\documentclass[11pt]{article}

\usepackage[usenames]{color}

\title{Approximate and pseudo-amenability of various classes of Banach algebras}
\author{Y. Choi, F. Ghahramani, Y. Zhang}
\date{February 12, 2009.\footnote{
\ackndiag\ A~revised version is to appear in {\it J.~Funct.~Anal}\/.}}

\usepackage{amsmath,amsfonts}
\numberwithin{equation}{section}

\usepackage{diagrams}
\newcommand{\ackndiag}{{Uses Paul Taylor's {\tt diagrams.sty} macros.}}
\newarrowmiddle{x}{\times}{\times}{\times}{\times}
\newarrow{Dash} {}{dash}{}{dash}>
\newarrow{Eq} ===== 
\newarrow{No} {}{dash}{x}{dash}>  
\newarrow{twoway} <---> 

\usepackage{enumerate}

\newcommand{\appctr}{approx\-imately con\-tract\-ible}
\newcommand{\appam}{approx\-imately amen\-able}
\newcommand{\appctry}{approx\-imate con\-tract\-ib\-ility}
\newcommand{\appamy}{approx\-imate amen\-ab\-ility}
\newcommand{\appdiag}{approx\-imate diag\-onal}

\newcommand{\al}{\alpha}
\newcommand{\gm}{\gamma}
\newcommand{\Gm}{\Gamma}
\newcommand{\kp}{\kappa}
\newcommand{\lm}{\lambda}
\newcommand{\Lm}{\Lambda}

\newcommand{\veps}{\varepsilon}
\newcommand{\Ind}{{\mathcal I}}

\newcommand{\indic}[1]{{{\bf 1}_{#1}}}
\newcommand{\Lmmax}{\Lambda_\vee}
\newcommand{\Imax}{\Ind_\vee}
\newcommand{\Alg}[1]{{\ell^1(#1)}}
\newcommand{\Fin}{\mathop{\rm FIN}}

\newcommand{\Cplx}{{\mathbb C}}
\newcommand{\Nat}{{\mathbb N}}

\newcommand{\Zahl}{{\mathbb Z}}

\newcommand{\defeq}{:=}
\newcommand{\dt}[1]{\textcolor{Bittersweet}{\sf#1}}
\newcommand{\st}{\;:\;}
\newcommand{\iso}{\cong}
\newcommand{\op}{{\mathop{\rm op}}}
\newcommand{\wtild}[1]{\widetilde{#1}}
\newcommand{\tp}{\otimes}
\newcommand{\ptp}{\widehat{\otimes}}
\newcommand{\wstar}{{\ensuremath{{\rm w}^*}}}

\newcommand{\clin}{\overline{\mathop{\rm lin}}}
\newcommand{\abs}[1]{{\left\vert#1\right\vert}}
\newcommand{\norm}[1]{{\Vert#1\Vert}}

\newcommand{\pair}[2]{\langle#1\/,\/#2\rangle}

\newcommand{\lp}[1]{\ell^{#1}}
\newcommand{\Lp}[1]{L^{#1}}
\newcommand{\lip}{{\mathop{\rm lip}\nolimits}} 
\newcommand{\Bdd}{{\mathcal B}}
\newcommand{\id}[1][]{{\sf 1}_{#1}}

\newcommand{\fu}[1]{{#1}^\#}
\newcommand{\Cmm}{\mathop{\rm comm.}}
\newcommand{\dfn}{\scriptstyle\mathop{\sf def}}

\newcommand{\bbI}{{\mathbb I}}
\newcommand{\bbJ}{{\mathbb J}}

\newcommand{\fg}[1]{{\mathbb F}_{#1}}
\newcommand{\sph}[1]{{\mathbb S}(#1)}
\newcommand{\sA}{{\sf A}}
\newcommand{\sB}{{\sf B}}

\newcommand{\sN}{{\sf N}}

\newcommand{\om}{\omega}

\newcommand{\gA}{g} 
\newcommand{\gB}{h} 

\newcommand{\LG}{{L^1(G)}} 
\newcommand{\SG}{{S^1(G)}}  

\newcommand{\redC}[1]{C_r^*(#1)}

\newcommand{\psf}[1]{PF_{#1}}  

\RequirePackage{amsthm}   

 \newcounter{pulse}[section]
\numberwithin{pulse}{section}  

\theoremstyle{plain}
\newtheorem{thm}[pulse]{\sc Theorem}
\newtheorem{propn}[pulse]{\sc Proposition}
\newtheorem{lemma}[pulse]{\sc Lemma}
\newtheorem{coroll}[pulse]{\sc Corollary}
\theoremstyle{definition}
\newtheorem{defn}[pulse]{\sc Definition}
\newtheorem*{notn}{\sc Notation}
\newtheorem{eg}[pulse]{\sc Example}
\newtheorem{egs}[pulse]{\sc Examples}
\theoremstyle{remark}

\newtheorem*{remstar}{\sc Remark}

\begin{document}
\maketitle
\begin{abstract}
We continue the investigation of notions of approximate amenability that were introduced in work of the second and third authors together with R.~J.~Loy. It is shown that every boundedly approximately contractible Banach algebra has a bounded approximate identity,
and that the Fourier algebra of the free group on two generators is not operator approx\-imately amenable.
Further examples are obtained of $\ell^1$-semigroup algebras which are approx\-imately amenable but not amenable;
 using these, we show that bounded approx\-imate contractibility need not imply sequential approx\-imate amenability.
Results are also given for Segal
algebras on locally compact groups, and 
algebras of $p$-pseudofunctions on discrete groups.

\end{abstract}

\begin{section}{Introduction}
In this article we continue the investigation of various notions of \appamy, initiated in \cite{GhL_genam1} and continued in papers by various authors: see, for instance, \cite{DLZ_St06,GhZ_pseudo,LashSam_SF05,Pir_appproj}.
Most of this paper is taken up with consideration of certain classes of Banach algebras, such as Fourier algebras, Segal subalgebras of $\LG$\/, certain $\lp{1}$-semigroup algebras and the algebras $PF_p(\Gm)$ where $\Gm$ is a discrete group, and the problem of determining when such algebras are \appam\ or pseudo-amenable (the definitions will be given below).

\subsection*{Overview of contents}
After establishing some background definitions and notation in Section~\ref{s:prelim}, we discuss some results for general Banach algebras.
If $A$ is a Banach algebra with an approximate diagonal, the forced unitization $\fu{A}$ need not possess an approximate diagonal (for example we may take $A$ to be $\ell^1$ with pointwise multiplication and apply \cite[Theorem~4.1]{DLZ_St06}). On the other hand, it follows from \cite[Proposition 3.2]{GhZ_pseudo} that if $A$ has a bounded approximate identity and an approximate diagonal, then $\fu{A}$ also has an approximate diagonal.
 We present a partial extension of this result to the case of multiplier-bounded {\appdiag}s: namely, we show that if $A$ has a central b.a.i.\ and a multiplier-bounded \appdiag, then so does~$\fu{A}$.

One outstanding basic question in this area is the following: does every \appam\ algebra have a bounded approximate identity? Although we are not able to resolve this here, we obtain some general results (Section \ref{s:general}) showing that slightly stronger notions of \appamy\ guarantee the existence of a bounded approximate identity.
As a consequence we are able to show
that several classes of Banach algebras, which might plausibly be pseudo-amenable, cannot be boundedly \appam: these include the Schatten classes ${\mathcal S}_p$ for $1\leq p < \infty$,
and all proper Segal algebras on locally compact groups.
A related argument shows that for any infinite compact metric space $X$\/, the Lipschitz algebra $\lip_\al(X)$\/, $0<\al\leq 1$\/, is not boundedly \appctr.

The last four sections are largely independent of each other and can be read interchangeably.
Section~\ref{s:AF2} resolves a question from \cite{GhSt_AG}, by showing that the Fourier algebra of $A(\fg{2})$ is not (operator) \appam. The proof uses direct manipulation of norm estimates, which rely on the ``rapid decay'' estimates known for $\fg{2}$ and $\fg{2}\times\fg{2}$\/. It then follows from known restriction theorems that whenever a locally compact group $G$ contains $\fg{2}$ as a closed subgroup, $A(G)$ is not (operator) \appam.

Section \ref{s:segal} collects some results on approximate notions of amenability for Segal algebras on locally compact groups. It is observed that Feichtinger's Segal algebra on an infinite compact abelian group is not \appam. We also show that if $\SG$ is pseudo-contractible for some Segal subalgebra $\SG\subseteq\LG$\/, then $G$ must be compact. It is also shown that whenever $G$ is a SIN-group, every Segal subalgebra $\SG\subseteq\LG$ is \emph{approximately permanently weakly amenable}\/: our proof uses the recent solution by Losert to the derivation problem for measure algebras.

Section \ref{s:l1ordered} is devoted to the $\lp{1}$-convolution algebras of totally ordered sets: when the sets are infinite, these algebras are never amenable. We show that these algebras are always boundedly \appctr\ (and in particular are boundedly \appam), but need not be sequentially \appam. This strengthens the observation in \cite{GhLZ_genam2} that the convolution algebra $\Alg{\Omega_\wedge}$ is boundedly \appctr\ but not sequentially \appctr. 

Finally, in Section \ref{s:gpalg} we consider the algebras of $p$-pseudofunctions on discrete groups. As a special case of the results in this section, we show that if $\Gm$ is a discrete non-amenable group then its reduced $C^*$-algebra is not \appam. This gives some evidence for the tentative conjecture that all \appam\ $C^*$-algebras are automatically amenable. Further evidence is provided by the fact that not only are the $C^*$-algebras $\Bdd(H)$ and $\prod_n M_n(\Cplx)$ not amenable, but they are not even \appam\ -- this observation appears to be due to Ozawa, see the remark after Definition $1.2$ in \cite{Oz_nonam}. 

\subsection*{Acknowledgements}
F. Ghahramani was supported by NSERC grant 36640-07. Y. Zhang was supported by NSERC grant 238949-05.

\end{section}

\begin{section}{Preliminaries}\label{s:prelim}
\subsection{Definitions and notation}

Throughout, if $A$ is a Banach algebra we shall write $\fu{A}$ for the forced unitization of $A$. The adjoined identity element will usually be denoted by $\id$ unless stated otherwise.  

We will frequently use $\pi$ to denote the bounded linear map $A\ptp A\to A$ that is specified by $\pi(a\tp b)=ab$ ($a,b\in A)$\/.

Recall that a Banach algebra $A$ is said to be \dt{\appam} if for every $A$-bimodule $X$ and every bounded derivation $D:A\to X^*$ there exists a net $(D_\alpha)$ of \emph{inner} derivations such that $\lim_\alpha D_\alpha(a)=D(a)$ for all $a \in A$. 

$A$ is said to be:
\begin{itemize}
\item[--] \dt{boundedly \appam} if the net $(D_\alpha)$ can always be taken to be bounded (in the usual norm of ${\cal L}(A,X^*)$);
\item[--] \dt{sequentially \appam} if we can choose the net $(D_\alpha)$ to be a sequence.
\end{itemize}
By the uniform boundedness principle (or a more direct Baire category argument) one sees that sequential \appamy\ implies bounded \appamy. The converse is not in general true, as will be shown in Section \ref{s:l1ordered} by combining Theorems \ref{t:total-ord_is_BAC} and \ref{t:total-ord_not_CAC}.

A Banach algebra $A$ is \dt{\appctr}\ if for every continuous derivation $D:A\to X$, where $X$ is a Banach $A$-bimodule, there exists a net $(D_i)$ of inner derivations such that $\lim_i D_i(a)=D(a)$ for all $a\in A$. The corresponding variants of \dt{bounded} and \dt{sequential}
 \appctry\ are defined in analogous fashion to the corresponding notions of \appamy.

\begin{remstar}
It is shown in \cite[Theorem 2.1]{GhLZ_genam2} that the concepts of \appctry\ and \appamy\ are in fact equivalent. However, this is not true for the corresponding sequential variants, and remains unknown (at present) for the bounded variants.
\end{remstar}

It has proved very useful in the classical theory of amenability to have characterizations in terms of virtual diagonals or {\appdiag}s. In much of this paper we shall work with {\appdiag}s rather than nets of derivations. To fix terminology we recall the following definition.

\begin{defn}\label{dfn:MBAD}
Let $A$ be a Banach algebra. An \dt{\appdiag} for $A$ is a net $(M_i)$ in $A\ptp A$ such that, for each $a\in A$\/,
\[ aM_i- M_ia\to 0 \quad\text{ and }\quad a\pi(M_i)\to a \,.\]
We say that the diagonal $(M_i)$ is \dt{multiplier-bounded} if there exists a constant $K$ such that  for all $a\in A$ and all $i$, each of
\begin{equation}\label{eq:MBAD-dfn}
\norm{aM_i- M_ia} \,,\;  
\norm{a\pi(M_i)-a} \text{ and } \norm{\pi(M_i)a-a}
\end{equation}
is bounded by $K\norm{a}$\/.
\end{defn}

The following equivalence is easily verified.

\begin{propn}
A Banach algebra $A$ is boundedly \appctr\ if and only if $\fu{A}$ has a multiplier-bounded approximate diagonal.
\end{propn}

We shall also make brief use of the notions of pseudo-amenability and pseudo-contractibility. For convenience we recall the relevant definitions from \cite{GhZ_pseudo}.

\begin{defn}
Let $A$ be a Banach algebra. We say that $A$ is \dt{pseudo-amenable} if it has an \appdiag, and \dt{pseudo-contractible} if it has an \appdiag\ $(M_i)$ which satisfies $aM_i=M_ia$ for all $a\in A$ and all $i$\/.
\end{defn}

\subsection{Basics}
\begin{propn}\label{p:unitize}
Let ${\mathcal S}$ be one of the following classes of Banach algebras: \appam, 
\appctr, sequentially \appam, sequentially \appctr, boundedly \appam, boundedly \appctr.

Let $A$ be a Banach algebra. Then $A\in{\mathcal S}$ if and only if $\fu{A}\in{\mathcal S}$.
\end{propn}
\begin{proof}
The case of \appamy\ is given by \cite[Proposition~2.4]{GhL_genam1}, and in fact the proof follows the same line for all the other cases. The key points are that (i)~every derivation from $A$ can be extended to a derivation from $\fu{A}$, such that the extended derivation is inner if and only if the original one was; (ii)~if $D$ is a derivation from $\fu{A}$ to an $A$-bimodule $X$, and $e$ denotes the identity of $\fu{A}$, then there is an inner derivation $D_1:\fu{A}\to X$ such that $(D-D_1)(e)=0$\/.
\end{proof}

\begin{remstar}
Note that the proofs of ``$A$ \appctr $\iff$ $\fu{A}$ \appctr'' and $A$ \appam $\iff$ $\fu{A}$ \appam'' do not rely on the fact that \appctry\ and \appamy\ are equivalent.
\end{remstar}

\begin{thm}\label{t:BAC_intrinsic}
Let $A$ be a boundedly \appctr\ Banach algebra. Then there exists a constant $C>0$ and nets $(M_i)$ in $A\ptp A$ and nets $(F_i)$, $(G_i)$ in $A$ such that
\begin{enumerate}[$(i)$]
\item $\pi(M_i)=F_i+G_i$;
\item\label{cnd:RAI}
	 $aF_i\to a$ for all $a\in A$;
\item\label{cnd:RMB}
	$\norm{aF_i}\leq C\norm{a}$ for all $a\in A$ and all $i$;
\item $G_ia\to a$ for all $a\in A$;
\item\label{cnd:LMB}
	$\norm{G_ia}\leq C\norm{a}$ for all $a\in A$ and all $i$;
\item\label{eq:mainest}
	for all $a\in A$ and all $i$, $\norm{aM_i-M_ia -a\tp G_i + F_i \tp a} \leq C\norm{a}$;
\item\label{eq:mainconv}
	for all $a\in A$,
	$aM_i-M_ia -a\tp G_i + F_i \tp a \to 0$.
\end{enumerate}
\end{thm}

For the sake of completeness we give the proof.
\begin{proof}
Regard $\fu{A}\ptp\fu{A}$ as an $A$-bimodule in the usual way. Let $K$ be the kernel of the product map $\fu{A}\ptp\fu{A}\to \fu{A}$ and let $D:A\to K$ be the derivation defined by $D(a)=a\tp\id-\id\tp a$.

Since $A$ is boundedly \appctr, there exists a net $(u_i)$ in $K$ such that
\[ C\defeq \sup_i \sup_{\norm{a}\leq 1} \norm{au_i-u_ia} < \infty \]
and $au_i-u_ia\to D(a)$ for all $a\in A$.

Identifying $\fu{A}\ptp\fu{A}$ with the direct sum $A\ptp A \oplus A\tp\id \oplus \id\tp A \oplus \Cplx\id\tp\id$, we may write each $u_i$ in the form
\[ u_i = (-M_i) \oplus (F_i\tp\id) \oplus (\id\tp G_i) \]
for some $M_i\in A\ptp A$ and some $F_i, G_i\in A$. We shall show that the nets $(M_i)$, $(F_i)$ and $(G_i)$ have the required properties.

First, note that since $u_i\in K$ for all $i$ we must have
\[ 0=\pi(u_i) = -\pi(M_i) +F_i +G_i \quad\text{ for all $i$}.\]

Next: since
\[ au_i - u_i a = (-aM_i + M_ia +a\tp G_i -F_i\tp a) \oplus aF_i\tp\id \oplus (-\id\tp G_i a) \]
and the left-hand side is bounded in norm by $C\norm{a}$ for all $a$, we must have
$\norm{aF_i}\leq C\norm{a}$, $\norm{G_ia}\leq C\norm{a}$ and
\[ \norm{aM_i-M_ia -a\tp G_i + F_i\tp a}\leq C\norm{a}\]
for all $i$ and all $a$. 

Finally, for each $a\in A$ we have
\[ \begin{aligned}
a\tp \id-\id\tp a & = D(a) \\
& = \lim_i au_i-u_i a \\
& = \lim_i \left( (-aM_i + M_ia +a\tp G_i -F_i\tp a) \oplus aF_i\tp\id \oplus (-\id\tp G_i a)\right)
\end{aligned} \]
and matching up terms we conclude that
\[ a=\lim_i aF_i = \lim_i G_ia \quad\text{ and }\quad 0 = \lim_i aM_i - M_ia -a\tp G_i +F_i\tp a \]
as required.
\end{proof}

\begin{remstar}
It follows from this that every boundedly \appctr\ Banach algebra has a multiplier-bounded right approx\-imate identity and a multiplier-bounded left approx\-imate identity. We shall use this later, in Section~\ref{s:general}.
\end{remstar}

Let $\kappa$ denote the canonical embedding of $A$ into $A^{**}$\/.
We have the following analogue of Theorem \ref{t:BAC_intrinsic}.
\begin{thm}\label{t:BAA_intrinsic}
Let $A$ be a boundedly \appam\ Banach algebra. Then there exists a constant $C>0$ and nets $(M_i)$ in $(A\ptp A)^{**}$ and nets $(F_i)$, $(G_i)$ in $A^{**}$ such that
\begin{enumerate}[$(i)$]
\item $\pi(M_i)=F_i+G_i$;
\item\label{cnd:bidRAI}
	 $aF_i\to \kappa(a)$ for all $a\in A$;
\item\label{cnd:bidRMB}
	$\norm{aF_i}\leq C\norm{a}$ for all $a\in A$ and all $i$;
\item $G_ia\to \kappa(a)$ for all $a\in A$;
\item\label{cnd:bidLMB}
	$\norm{G_ia}\leq C\norm{a}$ for all $a\in A$ and all $i$;
\item\label{eq:bidmainest}
	for all $a\in A$, $\sup_i\norm{aM_i-M_ia -a\tp G_i + F_i \tp a} \leq C\norm{a}$;
\item\label{eq:bidmainconv}
	for all $a\in A$, $aM_i-M_ia -a\tp G_i + F_i \tp a \to 0$.
\end{enumerate}
\end{thm}
We omit the proof: the argument exactly follows the one for Theorem~\ref{t:BAC_intrinsic}.

\subsection{Two lemmas using approximate diagonals}
We record some lemmas here which will be used later. Both are natural adaptations of routine arguments from the setting of \emph{amenable} Banach algebras.

\begin{lemma}\label{l:extension}
Let $B$ be a unital Banach algebra with identity element $\id$, $A\subseteq B$ a closed subalgebra that contains $\id$,
and suppose that there exists a tracial functional $\tau$ on $A$ such that $\tau(\id)=1$. If $A$ is pseudo-amenable, then there exists a net $(\psi_\al)$ in $B^*$ such that $\psi_\al(\id)\to 1$ and
\[ \sup_{b\in B} \abs{\psi_\al(ab-ba)} \to 0 \qquad\text{ for any $a\in A$.} \]
\end{lemma}
Note that by a trivial rescaling, the net $(\psi_\al)$ in the conclusion of our lemma can be chosen such that $\psi_\al(\id)=1$ for all $\al$\/. However, the formally weaker property $\psi_\al(\id)\to \id$ will suffice for our intended application.

\begin{proof}
Let $(u_\al)$ be an \appdiag\ for $A$\/: note that since $A$ has an identity element $\id$\/, $\pi(u_\al)\to\id$\/. For each $\al$ we may write $u_\al = \sum_i c^{\al}_i \tp d^{\al}_i$\/,
where $c^\al_i,d^\al_i\in A$ for all $i$ and
$\sum_i\norm{c^\al_i}\norm{d^\al_i}<\infty$\/.
Let $\wtild{\tau}\in B^*$ be any bounded extension of $\tau$ to a functional on~$B$, and define
\begin{equation*}
\psi_\al(S) = \wtild{\tau}\left(\sum_i d^\al_i S c^\al_i \right) \qquad(S\in B).
\end{equation*}
Then \emph{since $\tau$ is a trace on $A$},
\[ \psi_\al(\id) = \tau\left(\sum_i d^\al_i c^\al_i \right) =
 \tau\left(\sum_i c^\al_i d^\al_i \right) = \tau\pi(u_\al) \to \tau(\id) \]
 and by hypothesis $\tau(\id)=1$.

For fixed $b\in B$, define a functional $\phi_b\in (A\ptp A)^*$ by
$$ \phi_b (x\tp y) = \widetilde{\tau}(ybx) \qquad(x,y\in A).$$
By definition of the projective tensor norm, we have $\norm{\phi_b} \leq \norm{\widetilde{\tau}}\/\norm{b}$.

Now for each $a\in A$, the tracial property of $\tau$ gives
\[\phi_b(u_\al a) = \phi_b\left(\sum_i c^\al_i \tp d^\al_i a\right)
= \widetilde{\tau}\left(\sum_i d^\al_i ab c^\al_i \right)
= \psi_\al(ab) \]
and
\[\phi_b(au_\al) = \phi_b\left(\sum_i ac^\al_i \tp d^\al_i\right)
= \widetilde{\tau}\left(\sum_i d^\al_i ba c^\al_i \right)
= \psi_\al(ba)\/. \]
Therefore
\begin{align*}
 \sup_{b\in B, \norm{b}\leq 1} \abs{ \psi_\al(ab-ba)}
 & = 
 \sup_{b\in B, \norm{b}\leq 1} \abs{ \phi_b(au_\al-u_\al a)} \\
 & \leq \norm{\wtild{\tau}} \/ \norm{au_\al -u_\al a} \to 0
\end{align*}
for each $a \in A$\/.
Thus $(\psi_\al)$ has the required properties.
\end{proof}

Our second lemma will be needed for the proof of Theorem \ref{t:total-ord_not_CAC}.
It says, loosely, that the Gelfand transforms of an approximate diagonal must converge pointwise to the indicator function of the diagonal of the character space.

\begin{lemma}\label{l:Gelf-tr}
Let $A$ be a Banach algebra with non-empty character space $\Phi_A$, and suppose $A$ has a bounded approximate identity.

If $A$ is \appam, then there exists a net $(\Delta_\alpha)\in (A\ptp A)^{**}$ with the following properties:
\begin{enumerate}[$(i)$]
\item $\lim_\alpha\pair{\Delta_\alpha}{\varphi\tp\chi} = 0$ for any $\varphi,\chi\in\Phi_A$ with $\varphi\neq\chi$;
\item $\pair{\Delta_\alpha}{\varphi\tp\varphi}=1$ for all $\alpha$ and any $\varphi \in \Phi_A$.
\end{enumerate}
Moreover, if $A$ is sequentially \appam, we can take $\Delta_\alpha$ to be a sequence.
\end{lemma}

Before giving the proof we note that one could shorten the argument slightly by appealing to a \emph{modification of} the proof of \cite[Proposition 3.2]{GhZ_pseudo}. However, since that proposition does not deal explicitly with the \emph{sequential} case, which will be crucial in our intended application, we have chosen a more direct and only slightly longer argument.

\begin{proof}
We shall only prove the statement in the case where $A$ is sequentially \appam\ (the case where we merely assume $A$ to be \appam\ is completely analogous).

Thus, suppose $A$ has a bounded approximate identity and is sequentially \appam. Let $E$ be any weak$^*$-limit point in $A^{**}$ of the bounded approximate identity of $A$, so that $aE=Ea=\kp(a)\in A^{**}$, $\kp$ denoting the embedding of $A$ in its bidual.

Let $\pi:A\ptp A\to A$ be the product map and let $K=\ker\pi$; this is a sub-$A$-bimodule of $A\ptp A$. Define a bounded derivation $D:A\to K^{**}$ by $D(a)=a\tp E - E\tp a$ ($a\in A$\/).
Since $A$ is sequentially \appam\ there
 exists a sequence $(u_n)$ in $K^{**}$ such that
$a\cdot u_n- u_n\cdot a\to D(a)$ for all $a\in A$ and all $n$\/.

Define $\Delta_n = E\tp E - u_n \in (A\ptp A)^{**}$\/. We have
\[ \pair{\Delta_n}{\varphi\tp\varphi} =\pair{\varphi}{E}^2 -\pair{u_n}{\varphi\tp\varphi} = 1 - \pair{u_n}{\pi^*(\varphi)} = 1 - \pair{\pi^{**}(u_n)}{\varphi} = 1 \]
for all $n$. Moreover, if $\varphi$ and $\chi$ are distinct characters on $A$\/, then there exists $a\in A$ with $\varphi(a)\neq \chi(a)$. Then
\begin{align*}
 \pair{a\cdot u_n-u_n\cdot a}{\varphi\tp\chi} 
 & = \pair{u_n}{(\varphi\tp\chi)\cdot a} - \pair{u_n}{a\cdot(\varphi\tp\chi)} \\
 & = \pair{u_n}{\varphi(a)\varphi\tp\chi}- \pair{u_n}{\chi(a)\varphi\tp\chi} \\
 & = (\varphi(a)-\chi(a))\pair{u_n}{\varphi\tp\chi}\/,
\end{align*}
while
\[ \pair{D(a)}{\varphi\tp\chi} = \pair{a\tp E - E\tp a}{\varphi\tp\chi} =
 \varphi(a)-\chi(a) \]
so that, since $a\cdot u_n-u_n\cdot a \to D(a)$,
\[\varphi(a)-\chi(a) = (\varphi(a)-\chi(a))\lim_n\pair{u_n}{\varphi\tp\chi} .\]
Since $\varphi(a)-\chi(a)\neq 0$, this implies that
$1 = \lim_n\pair{u_n}{\varphi\tp\chi}$, and so
\[ \lim_n \pair{\Delta_n}{\varphi\tp\chi} = \pair{E\tp E}{\varphi\tp\chi}-\lim_n\pair{u_n}{\varphi\tp\chi} = 0\,, \]
as required.
\end{proof}

\end{section}
\begin{section}{General results}\label{s:general}
Recall that $A$ is \appctr\ if and only if $\fu{A}$ has an \appdiag\ (this is \cite[Proposition 2.6(a)]{GhL_genam1}).

We would like to have a better understanding of just when the presence of an \appdiag\ in $A$ guarantees an \appdiag\ in $\fu{A}$, and to obtain corresponding results for multiplier-bounded {\appdiag}s. Note that by combining the proof of $(ii)\implies(iii)$ in \cite[Proposition 3.2]{GhZ_pseudo} with $(3)\implies(1)$ of \cite[Theorem 2.1]{GhLZ_genam2}, we obtain the following result.

\begin{propn}
Let $A$ be a Banach algebra which has a bounded approximate identity and an \appdiag. Then $A$ is \appctr, and so $\fu{A}$ has an \appdiag.
\end{propn}

A natural hope is that the result just stated remains true if we replace `\appdiag' with `multiplier-bounded \appdiag'. We have been unable to verify this: the problem seems to be that while \appamy\ implies \appctry, it is not known if bounded \appamy\ implies bounded \appctry. The following result gives some partial answers.


\begin{propn}\label{p:PsA-BAC}
Let $A$ be a Banach algebra with a central b.a.i.\ $(e_\lm)_{\lm\in\Lm}$.
Suppose that $A$ has a multiplier-bounded \appdiag\ $(M_i)$\/. Then $\fu{A}$ has a multiplier-bounded \appdiag, and so $A$~is boundedly \appctr.
\end{propn}

\begin{proof}
Throughout we denote the adjoined  unit of $\fu{A}$ by $\id$\/, and the linearized product map $\fu{A}\ptp\fu{A}\to\fu{A}$ by $\pi$\/.
 We shall abuse notation and also use $\pi$ to denote the restricted map $A\ptp A\to A$\/.

We shall construct a net $(n_j)$ in $\fu{A}\ptp\fu{A}$ and a constant $K>0$ such that:
\begin{equation}\label{eq:Mike}
\norm{\pi(n_j)-\id}\leq K \text{ and } \norm{b\cdot n_j-n_j\cdot b} \leq K\norm{b} \text{ for all $b\in A$ and all $j$;}
\end{equation}
and
\begin{equation}\label{eq:Tires}
\lim_j \pi(n_j) =\id \text{ and }\lim_j (b\cdot n_j-n_j\cdot b) = 0\quad\text{ for all $b\in A$ and all $j$\/.}
\end{equation}
If these properties are satisfied, it is then straightforward to show that the net $(n_j)$ has the required properties in Definition~\ref{dfn:MBAD}.

By hypothesis there exist constants $C$ and $K_1$ such that $\norm{e_\lm}\leq C$ for all $\lm$ and such that, for all $a\in A$ and all $i$,
\begin{equation}\label{eq:oldMBAD1} 
\norm{a\pi(M_i)-a} \text{ and } \norm{a\cdot M_i-M_i\cdot a} \text{ are $\leq K_1\norm{a}$\/.}
\end{equation}
Moreover, for any $a\in A$, we have
\begin{equation}\label{eq:oldMBAD2} 
\lim_i a\pi(M_i)=a \text{ and } \lim_i a\cdot M_i-M_i\cdot a=0\/.
\end{equation}

To simplify the ensuing formulas slightly, we let $u_\lm\defeq 2e_\lm-e_\lm^2$ for each $\lm$\/: note that
$u_\lm + (\id-e_\lm)^2 = \id$ and $\norm{u_\lm}\leq 2C+C^2$\/,
for all $\lm$\/. We now set
\begin{equation}\label{eq:define}
m_{\lm,i} \defeq u_\lm \cdot M_i + (\id-e_\lm)\tp(\id-e_\lm)\/,
\end{equation}
so that
\begin{equation}\label{eq:auxil}
\pi(m_{\lm,i}) = u_\lm \pi(M_i) + \id - u_\lm\/.
\end{equation}

Let $\Ind$ and $\Lm$ be the index sets for the nets $(M_i)$ and $(e_\lm)$, respectively. We construct the required net $(n_j)$ using an iterated limit construction (see page 26 of \cite{Kelley_GT}). Our indexing directed set is defined to be $J=\Lm\times \prod_{\lm\in \Lm} \Ind$, equipped with the product ordering, and for each $j=(\lm, f)\in J$ we define $n_j=m_{\lm, f(\lm)}$. 

Fix $\lm$ and $i$. Using \eqref{eq:oldMBAD1} and \eqref{eq:auxil} gives
\[ \norm{\pi(m_{\lm,i})-\id} = \norm{u_\lm\pi(M_i)-u_\lm}\leq K_1\norm{u_\lm}\leq K_1(2C+C^2)\,.\]
Also, since each $e_\lm$ lies in the centre of $A$, we have for any $b\in A$
the identity
\begin{equation}\label{eq:newMBAD2}
b\cdot m_{\lm,i}-m_{\lm,i}\cdot b = u_\lm b\cdot M_i - u_\lm\cdot M_i\cdot b + (b-be_\lm)\tp (\id-e_\lm) - (\id-e_\lm)\tp(b-be_\lm)\,.
\end{equation}
Using \eqref{eq:oldMBAD1} again gives
\[ \norm{b\cdot m_{\lm,i}-m_{\lm,i}\cdot b} \leq \norm{u_\lm}\/\norm{b\cdot M_i-M_i\cdot b} + 2\norm{b}\/\norm{\id-e_\lm}^2 \leq (CK_1+2(1+C)^2)\norm{b} \]
for any $b\in A$. Since $\lm$ and $i$ were arbitrary, we have shown that \eqref{eq:Mike} holds with, say, $K=2(1+K_1)(1+C)^2$\/.

It remains to show that \eqref{eq:Tires} holds. 
Using \eqref{eq:oldMBAD2} and \eqref{eq:auxil} we have, for every $\lm$,
\[ \lim_i \pi(m_{\lm,i}) = \lim_i u_\lm\pi(M_i) +\id- u_\lm = \id \/;\]
hence, by \cite[Theorem~2.4]{Kelley_GT}, 
\[ \lim_j \pi(n_j)=\lim_\lm\lim_i \pi(m_{\lm,i}) =\id\/.\]
Using \eqref{eq:oldMBAD2} and \eqref{eq:newMBAD2} we have, for every $\lm$,
\[ \lim_i b\cdot m_{\lm,i}- m_{\lm,i}\cdot b = 
- (b-be_\lm)\tp (\id-e_\lm) + (\id-e_\lm)\tp(b-be_\lm)\,;\]
therefore, since $(e_\lm)$ is a bounded approximate identity for $A$, applying \cite[Theorem~2.4]{Kelley_GT} we obtain
\[ \lim_j (b\cdot n_j-n_j\cdot b_j) = \lim_\lm\lim_i  
 ( b\cdot m_{\lm,i}- m_{\lm,i}\cdot b ) = 0 \/.\]
Thus \eqref{eq:Tires} holds and our proof is complete.
\end{proof}

\begin{remstar}
The result is false if we do not require the central approximate identity in $A$ to be bounded: an example is given by $\lp{1}(\Nat)$ equipped with pointwise multiplication (\cite[Theorem 4.1]{DLZ_St06}).
\end{remstar}



It is still open whether an \appam\ Banach algebra must have a bounded approximate identity. If this were the case then one could extend many of the known hereditary properties of amenability to hold for \appamy.
All presently known examples of \appam\ Banach algebras have a bounded approximate identity.
In addition, all known examples of \appam\ Banach algebras are in fact \emph{boundedly \appctr}\/.

These last two observations are connected by the following result: every boundedly \appctr\ algebra has a bounded approximate identity (Corollary \ref{c:BAC-has-BAI} below).
We are able to prove a slightly stronger technical result, that allows us to rule out bounded \appamy\ for several classes of Banach algebras.
%


\begin{thm}\label{t:BAA+MBAI}
Suppose that the Banach algebra $A$ is boundedly \appam, and has both a multiplier-bounded left approximate identity and a multiplier bounded right approximate identity. Then $A$ has a bounded approximate identity.
\end{thm}

\begin{proof}
Let $(e_\al)$ and $(f_\beta)$ be, respectively, right and left multiplier-bounded approximate identities for $A$, so that there exists a constant $K>0$ such that
\[ \norm{ae_\al} \leq K \norm{a} \quad\text{ and }\quad \norm{f_\beta a}\leq K\norm{a}
  \quad\text{ for all $a\in A$ and all $\al,\beta$\/.} \]
From this we obtain the following estimates:
\begin{enumerate}[$(i)$]
\item $\norm{f_\beta\cdot m} \leq K\norm{m}$ and $\norm{m\cdot e_\al}\leq K\norm{m}$, for every $m\in A\ptp A$ and every $\al,\beta$\/;
\item $\norm{f_\beta\cdot T} \leq K\norm{T}$ and $\norm{T\cdot e_\al}\leq K\norm{T}$, for every $T\in (A\ptp A)^{**}$ and every $\al,\beta$\/.
\end{enumerate}
(The first pair of estimates follows easily from the definition of the projective tensor norm. The second pair follows from the first pair using Goldstine's theorem
 and the weak$^*$-continuity of the actions of $A$ on $(A\ptp A)^{**}$\/.)

Let $(F_i)$, $(G_i)$, $(M_i)$ and $C$ be the nets (respectively, the constant) satisfying $(ii)$--$(vi)$ of Theorem \ref{t:BAA_intrinsic}.

Suppose that the net $(f_\beta)$ is (norm) unbounded. We derive a contradiction as follows. For every $i$ and every $\beta$ we have
\[ \norm{f_\beta\cdot M_i - M_i\cdot f_\beta - f_\beta\tp G_i + F_i\tp f_\beta} \leq C \norm{f_\beta}\,, \]
and so by $(ii)$ above, we have
\begin{equation}\label{eq:Tim2}
 \norm{(f_\beta\cdot M_i - M_i\cdot f_\beta - f_\beta\tp G_i + F_i\tp f_\beta)e_\al} \leq KC \norm{f_\beta}\,, 
\end{equation}
for every $\al$\/, $\beta$ and~$i$\/.

Using the triangle inequality and the left multiplier-boundedness of the set $\{f_\beta\}$\/, from~\eqref{eq:Tim2} we have
\begin{equation}\label{eq:Daisy2}
\begin{aligned}
\norm{f_\beta}\/\norm{G_i\cdot e_\al}
 & \leq KC\norm{f_\beta} + \norm{f_\beta\cdot (M_i\cdot e_\al) } + \norm{M_i\cdot (f_\beta e_\al)} + \norm{F_i}\/\norm{f_\beta e_\al} \\
 & \leq KC\norm{f_\beta} + K\norm{M_i\cdot e_\al} + K\norm{M_i}\/\norm{e_\al} + K\norm{F_i}\/\norm{e_\al} \,,
\end{aligned}
\end{equation}
for every $\al$\/, $\beta$ and~$i$\/. Hence
\begin{equation}\label{eq:Brian2}
\norm{G_i\cdot e_\al} \leq KC+\frac{C}{\norm{f_\beta}}\left(
\norm{M_i\cdot e_\al} + \norm{M_i}\/\norm{e_\al} + \norm{F_i}\/\norm{e_\al}
\right) \,,
\end{equation}
for every $\al$\/, $\beta$ and~$i$\/.

For fixed $\al$ and $i$\/, combining \eqref{eq:Brian2} with our assumption that $\{ f_\beta\}$ is an unbounded set yields $\norm{G_i\cdot e_\al} \leq KC$\/. Taking limits with respect to~$i$\/, we then obtain $\norm{e_\al} \leq KC$ for each $\al$\/. But since $(e_\al)$ is a right approx\-imate identity and $(f_\gm)$ is a left-multiplier bounded set, we obtain
\[ \norm{f_\gm} = \lim_\al \norm{ f_\gm e_\al} \leq \lim_\al K \norm{e_\al} \leq K^2 C\/, \]
for all $\gm$\/. This contradicts our assumption that the net $(f_\beta)$ is unbounded.

A similar argument, with left and right interchanged, shows that the net $(e_\al)$ is also bounded; the existence of a bounded approximate identity is now standard.
\end{proof}

\begin{coroll}\label{c:BAC-has-BAI}
Let $A$ be boundedly \appctr. Then $A$ has a bounded approximate identity.
\end{coroll}

\begin{proof}
This is an immediate consequence of Theorem~\ref{t:BAA+MBAI}, since by Theorem \ref{t:BAC_intrinsic} every boundedly \appctr\ Banach algebra has a right and a left multiplier-bounded approximate identity.
\end{proof}

\begin{coroll}
Suppose that $A$ and $B$ are boundedly \appctr\ Banach algebras. Then the direct sum $A\oplus B$ is boundedly \appctr.
\end{coroll}

\begin{proof}
This follows from the proof of \cite[Proposition~2.7]{GhL_genam1} and Corollary~\ref{c:BAC-has-BAI}.
\end{proof}

The following result is similar to Theorem \ref{t:BAA+MBAI}, but seems not to imply it nor be implied by it.

\begin{propn}\label{p:BAA+MBset}
Suppose that the Banach algebra $A$ is boundedly \appam. Let $S$ be a subset of $A$ which is (left and right) multiplier bounded, i.e.~for some $K>0$\/, we have
$\norm{sa}\leq K\norm{a}$ and $\norm{as}\leq K\norm{a}$ for all $a\in A$, $s\in S$\/. Then $S$ is norm bounded.
\end{propn}

\begin{proof}[Proof sketch]
Arguing as at the start of the proof of Theorem \ref{t:BAA+MBAI}, we note that for every $s\in S$\/:
\begin{enumerate}[$(i)$]
\item $\norm{s\cdot m} \leq K\norm{m}$ and $\norm{m\cdot s}\leq K\norm{m}$, for every $m\in A\ptp A$\/;
\item $\norm{s\cdot T} \leq K\norm{m}$ and $\norm{T\cdot a}\leq K\norm{a}$, for every $T\in (A\ptp A)^{**}$\/.
\end{enumerate}

Suppose that $S$ is (norm) unbounded, so that there exists a sequence $(s_n)$ in $S$ with $\norm{s_n}\to\infty$\/. Then we may argue as in the proof of Theorem \ref{t:BAA+MBAI}, 
 replacing $e_\al$ with $s_m$ and $f_\beta$ with $s_n$ in Equations \eqref{eq:Tim2}--\eqref{eq:Brian2}, to show that the net $(s_m)$ is bounded, giving us a contradiction as before. Hence $S$ is (norm) bounded as claimed.
\end{proof}


\begin{egs}\label{eg:not_BAA}
The following algebras have multiplier-bounded approximate identities but have no bounded approximate identities.
\begin{enumerate}[$(a)$]
\item $c_0(\omega)$, the space of all sequences such that $\abs{a_n}\omega_n\to 0$, equipped with pointwise multiplication, where $\lim_n\omega_n=+\infty$.
\item $\lp{1}(\Nat_{\min},\omega)$, the weighted convolution algebra of the semilattice $\Nat_{\min}$, where $\lim_m\omega_n=\infty$.
\item The Schatten ideals ${\mathcal S}_p(H)$ ($H$ a Hilbert space) where $1\leq p < \infty$.
\item The Fourier algebras of \dt{weakly amenable}, non-amenable groups
 (see~\cite{deCann-Haa} for the def\-inition and examples).
\item Proper
symmetric
 Segal subalgebras (in the sense of Reiter \cite{Rei_Seg71}) of~$\Lp{1}(G)$\/.
\end{enumerate}
It therefore follows from Theorem~\ref{t:BAA+MBAI} that none of the above algebras can be boundedly \appam.
\end{egs}

\begin{remstar}
It has recently been shown (H.~G. Dales and R.~J. Loy, private communication) that the algebras of Example \ref{eg:not_BAA}(b) are not even \appam.
\end{remstar}

We can exploit Corollary~\ref{c:BAC-has-BAI} further to show that certain \emph{unital} Banach algebras are not boundedly \appctr.
\begin{coroll}
Let $X$ be an infinite, compact metric space and let $0 < \al \leq \frac{1}{2}$\/. Then the Lipschitz algebra $\lip_\al(X)$ is not boundedly \appctr.
\end{coroll}

\begin{proof}
Since $X$ is infinite and compact it contains a non-isolated point, $x_0$ say. Let $M=\{f\in\lip_\al(X)\st f(x_0)=0\}$\/: then $\fu{M}\iso\lip_\al(X)$\/. If $\lip_\al(X)$ were boundedly \appctr, then by Proposition~\ref{p:unitize} $M$ would also be boundedly \appctr, and hence by Corollary \ref{c:BAC-has-BAI} would have a bounded approximate identity. By Cohen's factorization theorem, this would imply that $M^2=M$\/, which is easily seen to be false by considering the function
$f: x \mapsto d(x,x_0)^{\beta}$ where $\al < \beta < 2\alpha$\/ (see also the remarks at the end of \cite{BCD_Lip}).
\end{proof}


\end{section}


\begin{section}{$A(\fg{2})$ is not \appam}\label{s:AF2}
Let $\fg{2}$ denote the free group on two generators. It was observed in \cite[Remark 3.4(b)]{GhSt_AG} that $A(\fg{2})$ is pseudo-amenable. In this section we shall show that it is not \appam; indeed, we prove the formally stronger result that $A(\fg{2})$ is not even operator \appam. Our techniques are based on direct estimates, exploiting the fact that the norm in $A(\fg{2}\times\fg{2})$ majorizes a certain weighted $\lp{2}$-norm. Some consequence for Fourier algebras of more general groups will be given at the end of the section. 

\subsection*{Background material}
We state the required definitions and basic properties in the setting of discrete groups, since we will eventually specialize to $\fg{2}$\/: some hold in greater generality, but we shall not discuss this here.

Let $\Gm$ be a discrete group and $C_{00}(\Gm)$ the space of compactly supported functions on~$\Gm$\/. The \dt{Fourier algebra} $A(\Gm)$ can be defined as the completion of $C_{00}(\Gm)$ with respect to the norm
\[ \norm{f}_{A(\Gm)} = \inf \{ \norm{\xi}_2\/\norm{\eta}_2 \st \xi,\eta\in \lp{2}(\Gm)\,;\, f=\xi*\eta \}
\quad(f\in C_{00}(\Gm)). \]

Let $\lm: \lp{1}(\Gm)\to B(\lp{2}(\Gm))$ denote the (faithful) left regular representation of $\lp{1}(\Gm)$ on $\lp{2}(\Gm)$\/. The WOT-closure of the image of $\lm$ is the \dt{von Neumann algebra of $\Gm$}, and will here be denoted by $VN(\Gm)$\/. We can identify $A(\Gm)$ with the predual of the group von Neumann algebra~$VN(\Gm)$\/: the pairing between the two satisfies
\[ \pair{\lm(T)}{f} = \sum_{g\in G} T(g)f(g) \qquad (f\in A(\Gm),\, T\in C_{00}(\Gm)\/), \]
from which the following is immediate.
\begin{lemma}\label{l:duality}
For every $f\in A(\Gm)$ and every $T\in C_{00}(\Gm)$ we have
\[ \abs{\sum_{x\in G} f(x) T(x) } \leq \norm{f}_{A(\Gm)} \, \norm{\lm(T)} \,.\]
\end{lemma}

The norm on $A(\Gm)$ is in general hard to describe, but
there are easy upper and lower bounds, which we record
 as a lemma for later reference.

\begin{lemma}\label{l:sandwich-bounds}
$\lp{2}(\Gm)\subseteq A(\Gm)$\/. Moreover,
\begin{equation}\label{eq:sandwich}
 \norm{f}_\infty \leq \norm{f}_{A(\Gm)} \leq \norm{f}_2 \qquad\text{ for all $f\in \lp{2}(\Gm)$\/.}
\end{equation}
\end{lemma}

For sake of completeness we give the proof.
\begin{proof}
If $f\in \lp{2}(\Gm)$ then $f=f*\delta_e$\/, where $e$ is the identity element of~$\Gm$\,: this proves the second inequality in \eqref{eq:sandwich}.

Suppose that $\xi,\eta\in\lp{2}(\Gm)$ satisfy $f=\xi*\eta$\/. Let $\check{\eta}$ be the element of $\lp{2}(\Gm)$ that is defined by $\check{\eta}(t)=\overline{\eta(t^{-1})}$\/, $t\in\Gm$\/. Then, given $x\in\Gm$\/, let $L_x: \lp{2}(\Gm)\to\lp{2}(\Gm)$ denote the operator of left translation by $x$\/, and note that
\[ f(x) =\xi*\eta(x) = \sum_{t\in\Gm} \xi(t)\eta(t^{-1}x)
  = \sum_{s\in\Gm} \xi(x^{-1}s) \eta(s^{-1})
 = \pair{L_x\xi}{\check{\eta}} \]
so that $\abs{f(x)} \leq \norm{L_x\xi}_2 \norm{\check{\eta}}_2 = \norm{\xi}_2\norm{\eta}_2$\/. Taking the supremum over all $x\in \Gm$ and the infimum over all pairs $(\xi,\eta)$ which `represent' $f$\/, we obtain the first inequality in \eqref{eq:sandwich}.
\end{proof}

\newcommand{\opptp}{{\widehat{\otimes}_{\rm op}}}

The following definition seems to have first appeared explicitly in \cite{GhSt_AG}.
\begin{defn}
Let $A$ be a quantized Banach algebra. We say that $A$ is \dt{operator \appam} if, for each quantized Banach $A$-bimodule~$X$, every completely bounded derivation $A\to X^*$ is approximately inner.
\end{defn}
Clearly, if $A$ is a quantized Banach algebra which happens to be \appam, then it is operator \appam.

The following is a `quantized' version of one direction of \cite[Proposition 3.3]{DLZ_St06}, specialized to the cases of interest.
\begin{propn}\label{p:qDLZ}
Let $\Gm$ be a discrete group and suppose that $A(\Gm)$ is operator \appam. Then
 for every finite set $S\subset A(\Gm)$ and every $\veps >0$\/,there exists $F\in c_{00}(\Gm\times \Gm)$ such that
\begin{itemize}
\item $\norm{a\cdot F - F\cdot a - a\tp\pi(F) + \pi(F)\tp a}_{A(\Gm\times\Gm)} \leq \veps$
\item $\norm{a-a\pi(F)}_{A(\Gm)}\leq\veps$
\end{itemize}
for every $a\in S$\/.
\end{propn}

For convenience we give a brief outline of how Proposition \ref{p:qDLZ} follows from existing results.

\begin{proof}[Proof sketch]
We use $\opptp$ to denote the operator projective tensor product of two operator spaces. Since $A(\Gm)$ is a quantized Banach algebra, $\fu{A(\Gm)}\opptp\fu{A(\Gm)}$ is a quantized Banach $A(\Gm)$-bimodule. Let $K$ be the kernel of the (surjective, completely bounded) product map $\fu{A(\Gm)}\opptp\fu{A(\Gm)}\to\fu{A(\Gm)}$\/: then $K$ and hence $K^*$ are also quantized Banach $A(\Gm)$-bimodules.

Let $D: A(\Gm) \to K^{**}$ be the completely bounded derivation defined by $D(a) =a\tp \id -\id\tp a$\, ($a\in A(\Gm)$\/). Since $K^{**}$ is the dual of a quantized Banach $A(\Gm)$-bimodule, by hypothesis $D$~is approximately inner. Therefore, by combining the proofs of \cite[Corollary~2.2]{GhL_genam1} and \cite[Proposition~2.1]{DLZ_St06}, we obtain the following: for any finite subset $S\subset A(\Gm)$ and any $\veps>0$\/, there exists $F\in A(\Gm)\tp A(\Gm)$ and $u,v\in A(\Gm)$ such that
\begin{itemize}
\item[(1)] $\norm{a\cdot F - F\cdot a + u\tp a -a \tp v}_{A(\Gm)\opptp A(\Gm)} < \veps$\/,
\item[(2)] $\norm{a-au}_{A(\Gm)} <\veps$ and $\norm{a-va}_{A(\Gm)}<\veps$\/.
\end{itemize}

By results of Effros and Ruan
\cite{ER90_OSAP},
 the operator projective tensor norm on $A(\Gm)\tp A(\Gm)$ coincides with its norm as a linear subspace of $A(\Gm\times\Gm)$\/. The rest of the proof now follows \cite[Propositions 2.3 and~3.3]{DLZ_St06} and we omit the details.
\end{proof}


\subsection*{Specializing to $\fg{2}$}

\begin{notn}
If $t\in\fg{2}$ then we denote by $\abs{t}$ the \dt{word length} of $t$\/, with the convention that the identity element has length~$0$\/. For each $n\in\Zahl_+$ let
\begin{equation}
 \sph{n} \defeq \{ t\in\fg{2} \st \abs{t}= n\}
\end{equation}
Elementary calculations show that for each $n\in\Nat$
\[ \abs{\sph{n}} = 4\cdot 3^{n-1} \;. \]
\end{notn}

The following Sobolev-type estimate, which we state without proof, is crucial for the argument to follow.
It is a special case of \cite[Theorem~1.1]{RRS_GAFA}, and as such really 
belongs to the province of geometric group theory,

 \begin{propn}\label{p:RRS_ineq}
Fix $m,n\in\Nat$\/. Let
 $T\in C_{00}(\fg{2}\times\fg{2})$ be supported on $\sph{m}\times\sph{n}$\/. Then
$\norm{\lm(T)} \leq (m+1)(n+1)\norm{T}_2$\;.
\end{propn}

\begin{coroll}
Let 
$F\in A(\fg{2}\times\fg{2})$\/. Then
\begin{equation}\label{eq:AF2F2_lowerbd}
\norm{F}_{A(\fg{2}\times\fg{2})} \geq \sup_{m,n\in\Zahl_+} \frac{1}{(m+1)(n+1)} \left( \sum_{x\in\sph{m}} \sum_{y\in\sph{n}} \abs{F(x,y)}^2\right)^{1/2} \;.
\end{equation}
\end{coroll}

\begin{proof}
This is a routine deduction from Proposition~\ref{p:RRS_ineq}, using duality.
Let $(m,n)\in\Zahl_+^2$ and let $T_{m,n}\in C_{00}(\fg{2}\times\fg{2})$ be defined by
\[ T_{m,n}(x,y) = \left\{ \begin{aligned}
(m+1)^{-1}(n+1)^{-1}\overline{F(x,y)} & \quad\text{ if $x\in\sph{m}$ and $y\in\sph{n}$} \\
0 & \quad\text{ otherwise.}
 \end{aligned} \right. \]
Then by Proposition \ref{p:RRS_ineq},
\[ \norm{\lm(T_{m,n})}\leq \left(\sum_{(x,y)\in\sph{m}\times\sph{n}} \abs{F(x,y)}^2\right)^{1/2} \,, \]
and since
\[ \abs{\sum_{(x,y)\in\fg{2}\times\fg{2}} F(x,y)T_{m,n}(x,y)} = \frac{1}{(m+1)(n+1)}\sum_{(x,y)\in\sph{m}\times\sph{n}} \abs{F(x,y)}^2\;, \]
applying Lemma \ref{l:duality} completes the proof.
\end{proof}

\subsection*{Proof that $A(\fg{2})$ is not operator \appam}
We start with some notation.
In view of the lower bound
\eqref{eq:AF2F2_lowerbd}, we introduce the following norm
on $c_{00}(\fg{2}\times\fg{2})$\/: given $H\in c_{00}(\fg{2}\times\fg{2})$\/, let
\[ \norm{H}_{\om\times\om} = \sup_{m, n\geq 0} \frac{1}{(m+1)(n+1)} \left(\sum_{x\in\sph{m}}\sum_{y\in\sph{n}}\abs{H(x,y)}^2\right)^{1/2}\;. \]

For each $n\in\Nat$, we fix a partition of $\sph{n}$ into two disjoint subsets $\sA(n)$ and $\sB(n)$ of equal cardinality, so that
\[ \abs{\sA(n)} = \abs{\sB(n)} = \frac{1}{2}\abs{\sph{n}}\,. \]
We also fix a sequence $(\gm_n)_{n\geq 1}$ of strictly positive reals, such that
\begin{equation}\label{eq:test_seq}
\sum_{n\geq 1} \gm_n^2 \abs{\sph{n}} < \infty
\end{equation}
(the $\gm_n$ will be chosen later with appropriate hindsight).
Now define elements $a$ and $b$ of $\lp{2}(\fg{2})$ by
\begin{equation}
a \defeq \sum_{m\geq 1} \gm_m \indic{\sA(m)}
 \quad,\quad
b \defeq \sum_{n\geq 1} \gm_n \indic{\sB(n)} \,.
\end{equation}
Finally, let $\veps_0>0$\/, $\veps_1>0$\/.

\medskip

{\bf Suppose that $A(\fg{2})$ is operator \appam.} Applying Proposition \ref{p:qDLZ} with $S=\{a,b\}$ and using the lower bounds from \eqref{eq:sandwich} and~\eqref{eq:AF2F2_lowerbd}, we obtain
 $F\in c_{00}(\fg{2}\times\fg{2})$ such that, if we set $u=\pi(F)\in c_{00}(\fg{2})$\/:
\begin{equation}\label{eq:start_AI}
\norm{a-ua}_\infty\leq \veps
\quad\text{ and }\quad \norm{b- ub}_{\infty} \leq \veps\,;
\end{equation}
\begin{subequations}
\begin{equation}\label{eq:start_a}
\norm{a\cdot F - F\cdot a -a\tp u+ u\tp a}_{\om\times\om} \leq \veps\,;
\end{equation}
\begin{equation}\label{eq:start_b}
\norm{b\cdot F - F\cdot b -b\tp u+ u\tp b}_{\om\times\om} \leq \veps\,.
\end{equation}
\end{subequations}

For the moment we shall ignore the relations \eqref{eq:start_AI}, and work exclusively with the information given by \eqref{eq:start_a} and \eqref{eq:start_b}. Our task will be simplified by the fact that we have chosen the functions $a$ and $b$ to have `large' yet disjoint supports (this theme, if not the actual calculations, is inspired by the proof of \cite[Theorem~4.1]{DLZ_St06}).

\begin{remstar}
 It is worth noting that, since $a$ and $b$ are fixed in advance of $F$\/,
 both \eqref{eq:start_a} and \eqref{eq:start_b}
 can always be satisfied by taking $F$ to be of the form $c\indic{W\times W}$ for some $c\in\Cplx$ and some suitably large, finite subset $W\subset \fg{2}$\/; hence we will need to use \eqref{eq:start_AI} at some point if we are to obtain the required contradiction.
\end{remstar}

Our task would be simplified if we furthermore assume that $F$ is constant on sets of the form $\sph{m}\times\sph{n}$\/, and indeed the calculations that follow are motivated by this special case.
The key step is contained in the following proposition.

\begin{propn}\label{p:almost_constant}
For each $k\in\Nat$ let
\[ \gA(k)\defeq \frac{1}{\abs{\sA(k)}} \sum_{p\in\sA(k)} u(p)
\quad\text{ and }\quad
 \gB(k)\defeq \frac{1}{\abs{\sB(k)}} \sum_{q\in\sB(k)} u(q) \]
Then for every $m,n\in\Nat$ we have
\begin{equation}\label{eq:master}
\abs{\sA(m)}^{1/2}\/\abs{\sB(n)}^{1/2}\/ \abs{ \gA(m)-\gB(n)} \leq (m+1)(n+1)(\gm_m^{-1}+\gm_n^{-1}) \veps_1
\end{equation}
and
\begin{subequations}
\begin{equation}\label{eq:small_var_a}
\abs{\sB(n)}^{1/2} \left( \sum_{x\in\sA(m)} \abs{u(x)-\gA(m)}^2 \right)^{1/2}
\leq (m+1)(n+1)(\gm_m^{-1}+\gm_n^{-1}) \veps_1
\end{equation}
\begin{equation}\label{eq:small_var_b}
\abs{\sA(m)}^{1/2} \left( \sum_{y\in\sB(n)} \abs{u(y)-\gB(n)}^2 \right)^{1/2}
\leq (m+1)(n+1)(\gm_m^{-1}+\gm_n^{-1}) \veps_1
\end{equation}
\end{subequations}
\end{propn}

For our proof we need a technical lemma, whose essential content is well-known, but is stated here for convenience.
\begin{lemma}\label{l:decouple}
Let $\bbI$, $\bbJ$ be finite index sets and let $c_i, d_j\in\Cplx$ for all $i\in\bbI$ and all $j\in\bbJ$\/. Let
\[ \mu_c\defeq \frac{1}{\abs{\bbI}} \sum_i c_i \quad\text{ and }\quad
   \mu_d\defeq \frac{1}{\abs{\bbJ}} \sum_j d_j \,.\] 
Then
\[ \frac{1}{\abs{\bbI}}\frac{1}{\abs{\bbJ}} \sum_{i\in\bbI}\sum_{j\in\bbJ} \abs{c_i-d_j}^2 = \abs{\mu_c-\mu_d}^2 + \frac{1}{\abs{\bbI}} \sum_{i\in\bbI} \abs{c_i -\mu_c}^2
 + \frac{1}{\abs{\bbJ}}\sum_{j\in\bbJ} \abs{d_j-\mu_d}^2 \,.\]
\end{lemma}

\begin{proof}[Proof sketch]
\newcommand{\Exp}{{\mathbb E}}
\newcommand{\Var}{{\mathop{\sf Var}}}
One can prove this by direct calculation. Alternatively, we can use the language of probability theory, as follows. If $X$ and $Y$ are independent complex-valued random variables defined on a common finite probability space (in this case, $\bbI\times\bbJ$) then 
\[ \begin{aligned}
\Exp\abs{X-Y}^2
 & = \Exp (X-Y)\overline{(X-Y)} \\
 & = \Exp X\overline{X} - (\Exp X)(\overline{\Exp Y}) - (\overline{\Exp X})(\Exp Y) + \Exp Y\overline{Y} \\
 & = (\Exp X - \Exp Y)^2 + \Exp \abs{X}^2 - \abs{\Exp X}^2 + \Exp \abs{Y}^2 - \abs{\Exp Y}^2
\end{aligned}\]
We also have
\[ \Exp \abs{ X- \Exp X}^2 = \Exp (X- \Exp X)(\overline{X- \Exp X}) = \Exp \abs{X}^2 - \abs{\Exp X}^2 \,,\]
and similarly $\Exp \abs{Y- \Exp Y}^2 =\Exp\abs{Y}^2 - \abs{\Exp Y}^2$\/. Taking $X$ to be the random variable $(i,j)\mapsto c_i$ and $Y$ to be the random variable $(i,j) \mapsto d_j$\/, the proof is complete.
\end{proof}

\begin{proof}[Proof of Proposition \ref{p:almost_constant}]
Equation \eqref{eq:start_a} implies that
\begin{equation*}
\begin{aligned}
  (m+1)^2(n+1)^2 \veps_1^2 
 & \geq \sum_{x\in\sph{m}}\sum_{y\in\sph{n}}
	\abs{ (a(x)-a(y))F(x,y) - a(x)u(y) + a(y)u(x)}^2 \\
 & \geq
	\sum_{x\in\sA(m),\,y\in\sB(n)}	\abs{ (a(x)-a(y))F(x,y) - a(x)u(y) + a(y)u(x)}^2 \\
 & = \sum_{x\in\sA(m)}\sum_{y\in\sB(n)} \abs{a(x)F(x,y)-a(x)u(y)}^2 \\
 & = \gm_m^2 \sum_{x\in\sA(m)}\sum_{y\in\sB(n)} \abs{F(x,y)-u(y)}^2 \;.
\end{aligned}
\end{equation*}
Therefore
\begin{subequations}
\begin{equation}\label{eq:gen_a}
 (m+1)(n+1)\veps_1 \geq \gm_m \left(\sum_{x\in\sA(m)}\sum_{y\in\sB(n)} \abs{F(x,y)-u(y)}^2\right)^{1/2} \,.
\end{equation}
Similarly, using Equation \eqref{eq:start_b} instead of \eqref{eq:start_a}, we have
\begin{equation}\label{eq:gen_b}
\begin{aligned}
 (m+1)(n+1)\veps_1
 & \geq \left( \sum_{x\in\sA(m)} \sum_{y\in\sB(n)} \abs{ -b(y)F(x,y) + b(y)u(x)}^2\right)^{1/2} \\
 & = \gm_n \left(\sum_{x\in\sA(m)} \sum_{y\in\sB(n)}
	\abs{F(x,y)-u(x)}^2\right)^{1/2}\;.
\end{aligned}
\end{equation}
\end{subequations}
Hence, by using the triangle inequality for the $2$-norm on $\lp{2}( \sA(m)\times\sB(n)\/)$\/, we see that \eqref{eq:gen_a} and \eqref{eq:gen_b} together imply
\begin{equation}\label{eq:BARRY}
(m+1)(n+1) (\gm_m^{-1}+\gm_n^{-1}) \veps_1
 \geq \left(
\sum_{x\in\sA(m)}\sum_{y\in\sB(n)} \abs{u(x)-u(y)}^2
\right)^{1/2}\;. 
\end{equation}
The desired estimates \eqref{eq:master}, \eqref{eq:small_var_a} and \eqref{eq:small_var_b} now follow by applying Lemma \ref{l:decouple} to~\eqref{eq:BARRY}.
\end{proof}

\bigskip

We now show that by fixing our sequence $(\gm_n)$ appropriately, we can force $g$ to be ``slowly varying at infinity'', and play this off against the fact that $g$ has finite support (since $u$ does). For each~$n$\/, take $\gm_n\defeq n^{-1} \abs{\sph{n}}^{-1/2}$
(this certainly satisfies the condition in~\eqref{eq:test_seq}\/).
If we substitute this into the estimate~\eqref{eq:master} and take $m=k$\/, $n=k+1$ we get
\[ \begin{aligned}
  & \abs{\sA(k)}^{1/2}\/\abs{\sB(k+1)}^{1/2}\/ \abs{\gA(k)-\gB(k+1)} \\
  & \leq \veps (k+1)(k+2) \left( k\abs{\sph{k}}^{1/2} + (k+1)\abs{\sph{k+1}}^{1/2} \right) \\
  & \leq \veps (k+2)^3 (\abs{\sph{k}}^{1/2}+\abs{\sph{k+1}}^{1/2}) \;;
\end{aligned} \]
and since
$\abs{\sA(n))} = \abs{\sB(n)} = \frac{1}{2}\abs{\sph{n}} = 2 \cdot 3^{n-1}$\,,
we find that
\[  \abs{\gA(k)-\gB(k+1)}
 \leq \veps(k+2)^3\, \frac{ 2\cdot 3^{(k-1)/2} + 2\cdot 3^{k/2} }{ 2\cdot 3^{(k-1)/2} \cdot 3^{k/2}}  = (1+\sqrt{3})\,\veps (k+2)^3\, 3^{-k/2} \;. \]
On the other hand, taking $m=n=k+1$ in \eqref{eq:master}, an exactly similar argument gives
\[ \abs{\gA(k+1)-\gB(k+1)} \leq 2\,\veps(k+2)^3\, 3^{-k/2} \,,\]
and we thus obtain the estimate
\begin{equation}\label{eq:geometric-tail_A}
 \abs{\gA(k+1)-\gA(k)} \leq 5\veps (k+2)^3\, 3^{-k/2}\;.
\end{equation}

\medskip
By the comparison test, the infinite sum $\sum_{k\geq 1} (k+2)^3 3^{-k/2}$ converges, with value $M$ say.
Moreover, since $u$ has finite support, there exists $N\geq 2$ such that $\gA(j)=\gB(j)=0$ for all $j\geq N$\/. Hence, using \eqref{eq:geometric-tail_A}, we get
\begin{equation}\label{eq:small_A}
\begin{aligned}
\abs{\gA(1)}
  = \abs{\gA(N)-\gA(1)} 
 & \leq \sum_{k=0}^{N-1} \abs{\gA(k+1) -\gA(k)} \\
 & \leq 5\veps \sum_{k=0}^{N-1} (k+2)^3 3^{-k/2}  
  < 5M\veps\;.
\end{aligned}
\end{equation}

\medskip
Now observe that, by \eqref{eq:start_AI},
\[ \veps
 \geq \norm{a-au}_\infty \geq \max_{x\in\sA(1)} \abs{a(x)-a(x)u(x)} = \frac{1}{2}\max_{x\in\sA(1)}\abs{1-u(x)}\,. 
 \]
Let $x,y$ be the two elements of $\sA(1)$\/. Then the estimate just given implies that
\[ 2\veps \geq \frac{1}{2}\left(\abs{1-u(x)}+\abs{1-u(y)}\right) \geq \frac{1}{2}\abs{2-u(x)-u(y)} = \abs{1-g(1)}\,,\]
and combining this with \eqref{eq:small_A}\/, 
we finally arrive at
\[ 1 \leq \abs{1-g(1)}+\abs{g(1)} \leq 2\veps+5 M\veps \,.\]
As $M$ is, by definition, \emph{independent of $\veps$}\/, we obtain a contradiction by taking $\veps$ to be sufficiently small.
Hence our assumption that $A(\fg{2})$ is operator \appam\ must be false, and the proof is complete.
\hfill$\Box$

\medskip

\begin{coroll}\label{c:cheapo}
Let $G$ be a locally compact group, into which $\fg{2}$ embeds as a closed subgroup\/. Then $A(G)$ is not operator \appam.
\end{coroll}

\begin{proof}
The hypothesis on $G$ ensures that the restriction homomorphism $A(G)\to A(\fg{2})$ is completely bounded and surjective (\cite[Theorem 1a]{Herz_AIF73}), hence a completely bounded quotient homo\-morphism. If $A(G)$ were operator \appam, then $A(\fg{2})$ would be also, since operator \appamy\ is inherited by completely bounded quotients. This gives a contradiction.
\end{proof}

\begin{remstar}
Let $G$ be a discrete group. If $G$ is amenable, then $A(G)$ has a bounded approximate identity (Leptin's theorem), and so is \appam\ by the arguments of \cite{GhSt_AG}.
\end{remstar}

Note also that there are discrete groups which are non-amenable yet contain no copy of $\fg{2}$\/: Ol'shanskii's groups, or Burnside groups of sufficiently large rank and exponent. So any attempt to prove that \appamy\ of $A(G)$ implies amenability of $G$ must use different, or additional, methods.
\end{section}


\begin{section}{Results for Segal algebras}\label{s:segal}
Following on from Example \ref{eg:not_BAA}$(e)$ above, we give some results on other notions of \appamy\ in the setting of Segal algebras.

Let $G$ be a locally compact group with a left-invariant Haar measure~$\lm$.
Throughout this section, $\SG$ denotes a Segal subalgebra of $\LG$ (in the sense of Rei\-ter \cite{Rei_Seg71}). We have already seen that $\SG$ is boundedly \appam\ if and only if it is equal to the whole of $\LG$ and $G$ is amenable. For Feich\-tinger's Segal algebra (see \cite{ReiSteg} for the definition) on a compact abelian group we easily obtain the following:

\begin{propn}
The Feich\-tinger algebra on an infinite compact abelian group is not \appam.
\end{propn}
\begin{proof}
When $G$ is compact and abelian, the Feich\-tinger algebra on $G$ is
\[  S_0(G) = \{f= \sum_{\gamma\in \hat G}c_\gamma\chi_\gamma \, :\; \norm{f} = \sum \abs{c_\gamma} < \infty \}\, ,  \]
where $\chi_\gamma$ is the character of $G$ associated with $\gamma\in \hat G$. Hence,
\[  S_0(G) \cong \ell^1(\hat G)  , \]
where the right-hand side is equipped with the pointwise product. But $\ell^1(S)$ is not \appam\ if $S$ is an infinite set, due to \cite[Theorem~4.1]{DLZ_St06}. So $S_0(G)$ is not \appam.
\end{proof}

It has been shown in \cite{GhZ_pseudo} that a Segal algebra on a compact group is pseudo-contractible. The converse is also true and is a consequence of the next proposition. 

\begin{propn}\label{p:central}
If there is $\sN\in \SG\ptp \SG$ such that $\pi(\sN)\neq 0$ and $f\cdot \sN = \sN\cdot f$ for $f\in \SG$, then $G$ is compact.
\end{propn}

\begin{proof} 

Let $\iota$: $\SG \to \LG$ be the inclusion injection. Then the following diagram commutes:
\[ \begin{diagram}
\SG\ptp\SG & \rTo^{\iota\tp\iota} & \LG\ptp \LG \\
\dTo^\pi & & \dTo_\pi \\
S (G) & \rTo_{\iota} & \LG
\end{diagram} \]
Let $\sN\in \SG\ptp \SG$ be such that $\pi(\sN)\neq 0$ and $f\cdot\sN = \sN\cdot f$ for $f\in \SG$. Let $M=\iota \otimes \iota(\sN)\in \LG\ptp \LG$\/. We have $f\cdot M = M\cdot f$ for all $f\in \LG$\/, and therefore $\mu\cdot M = M\cdot \mu$ ($\mu\in M(G)$). In particular, $M = \delta_{x^{-1}}\cdot M\cdot \delta_x$ ($x\in G$). Let $K$ be any compact subset of $G\times G$. If we regard $M$ as a function in $L^1(G\times G)$\/, then
\begin{align*}
  \int_K |M(s, t)|ds dt &= \int_K |\delta_{x^{-1}}\cdot M\cdot \delta_x (s, t)|ds dt\\ &= \int_{(x,e)K(e,x^{-1})}\Delta(x) |M(s, t)|ds dt \, ,
\end{align*}
where $(x,e)K(e,x^{-1})$ denotes the set $\{(xs, tx^{-1}): \; (s, t)\in K\}$. Given $\veps>0$, let $R\subset G\times G$ be a compact set such that
\[  \int_{G\times G \backslash R}|M(s, t)|ds dt < \veps \, . \]
If $G$ is not compact, then there is $x\in G$ such that $(x,e)K(e,x^{-1})\subset G\times G \backslash R$ and $\Delta(x) \leq 1$. So
\[  \int_{(x,e)K(e,x^{-1})}\Delta(x) |M(s, t)|ds dt < \veps \, .  \]
We then have $\int_K |M(s,t)|ds dt < \veps$. This implies that $\int_K |M(s,t)|ds dt = 0$ for all compact $K\in G\times G$. Therefore $M = 0$ in $\LG\ptp \LG$. On the other hand, $\pi(\sN) \neq 0$ in $\SG$ and hence $\pi(M)=\iota\pi(\sN) \neq 0$ in $\LG$, a contradiction. Thus, $G$ must be compact.
\end{proof}

\begin{remstar}
Proposition \ref{p:central} holds with $\SG$ replaced by any Banach algebra $B$ which admits a continuous injective homomorphism $B\to \LG$\/ whose range is dense. Therefore, if such a $B$ exists and is pseudo-contractible, $G$ must be compact.
\end{remstar}

Combining Proposition~\ref{p:central} and \cite[Theorem~4.5]{GhZ_pseudo}, we then have a characterization of a compact group.

\begin{thm} The following are equivalent for a locally compact group $G$.
\begin{enumerate}[$(i)$]
\item The group $G$ is compact;
\item there is a Segal algebra on $G$ which is pseudo-contractible;
\item all Segal algebras on $G$ are pseudo-contractible.
\end{enumerate}
\end{thm}

(A different proof of the part ``$(ii) \implies (i)$'' can be seen in \cite{S-S-S}.) 

It is natural to ask for an analogous characterization of amenability of $G$ in terms of \appamy\ or pseudo-amenability of Segal algebras on $G$\/. First we recall some material from the theory of abstract Segal algebras.

A dense left ideal $B$ of a Banach algebra $(A, \norm{\cdot}_A)$ is called an \dt{abstract Segal algebra in $A$}, or simply a \dt{Segal algebra in $A$}, with respect to some norm $\norm{\cdot}_B$ if it is a Banach algebra with respect to the norm $\norm{\cdot}_B$ and if $\norm{b}_A \leq \norm{b}_B$ ($b\in B$) \cite{Burnham1,Lein_Seg73}.  It was shown in \cite{Lein_Seg73} that if $B$ is a Segal algebra in $A$, then the mapping $J\mapsto \overline{J}^A$ is a bijection from the set of all closed right (two-sided) ideals in $B$ onto the set of all closed right (two-sided) ideals in $A$ and the inverse mapping is $I\mapsto I\cap B$, where for a set $J\subset B$ the notation $\overline{J}^A$ stands for the closure of $J$ in $A$. The same machinery as in \cite[Proposition~2.7]{Lein_Seg73} yields the following:

\begin{propn}\label{leinert}
Let $B$ be an abstract Segal algebra in a Banach algebra $A$, let $J$ be a closed ideal of $B$ and let $I = \overline{J}^A$.
\begin{enumerate}[$(i)$]
\item Suppose that $A$ and $J$ both have right approximate identities. The $I$ has a right approximate identity.
\item 
Suppose that $B$ and $I$ both have right approximate identities. Then~$J$ has a right approximate identity.
\end{enumerate}
\end{propn}

Since we only need part $(i)$ of Proposition \ref{leinert}, we shall give an independent proof of~$(i)$, which is more direct in the sense that it avoids the duality machinery of \cite{Lein_Seg73}.
\begin{proof}
Let $F\subset I$ be a finite subset, and let $\veps >0$. It suffices to find $s\in I$ such that $\max_{y\in F}\norm{ys-y}_A\leq\veps$.

Since $A$ has a right approx\-imate identity, there exists $u\in A$ with $\max_{y\in F}\norm{yu-u}_A < \veps/2$. Therefore, since $B$ is dense in $A$, there exists $u'\in B$ such that  $\max_{y\in F}\norm{yu'-y}_A < \veps/2$.

Since $I$ is a right ideal and $B$ is a left ideal in $A$, $yu'\in I\cap B=J$ for every $y\in F$. Therefore, as $J$ has a right approximate identity, there exists $w\in J$ such that
\[ \norm{yu'w-yu'}_B \leq \veps/2 \qquad\text{for all $y\in F$\/.} \]

Let $s=u'w$: since $w\in J\subseteq I$ and $I$ is also a \emph{left} ideal, we have $s\in I$; and for every $y\in F$,
\begin{align*}
 \norm{ys-y}_A
 & \leq \norm{yu'w-yu'}_A + \norm{yu'-y}_A \\
 & \leq \norm{yu'w-yu'}_B + \norm{yu'-y}_A  < \veps/2 + \veps/2 = \veps
\end{align*}
as required.
\end{proof}

\begin{remstar}
Part $(ii)$ of Proposition~\ref{leinert} can also be proved by direct $\veps$-$\delta$ arguments, similar to the ones just given; again, the duality machinery from \cite{Lein_Seg73} can be bypassed.
\end{remstar}

\begin{thm}\label{G amen}
If $\SG$ is \appam\ or pseudo-amenable then $G$ is an amenable group.
\end{thm}

\begin{proof}
Let $I_0=\{ f\in \LG \st \int_G f(x)\,dx=0\}$ be the augmentation ideal in~$\LG$, and let $J = I_0\cap \SG$. Then $J$ is a codim\-ension-$1$ two-sided closed ideal in $\SG$. If $\SG$ is \appam\ or pseudo-amenable, then $J$ has a right approximate identity by \cite[Corollary~2.4]{GhL_genam1} and \cite[Proposition~2.5]{GhZ_pseudo}.
By Proposition~\ref{leinert}$(i)$, $I_0$ must also have a right approximate identity. This implies that $G$ is amenable due to \cite[Theorem 5.2]{GAW_appid}. 
\end{proof}
\begin{remstar}
We do not know whether there is a Segal algebra $\SG$ that is \appam\ and that is not identical with $\LG$. Whether a Segal algebra on an amenable group is always pseudo-amenable is also open. Partial results can be found in \cite{GhZ_pseudo}.
\end{remstar}

We now turn to results that do not depend on amenability or compactness of~$G$\/. 
While $\LG$ is weakly amenable for every locally compact group~$G$\/, the same need not be true for Segal algebras unless $G$ is abelian. Following on from results in \cite{GhLau02, GhLau05} on approximate weak amenability of Segal algebras, we now look at approximate permanent weak amenability.


Recall from \cite{DGG_nWA}
that a Banach algebra $A$ is said to be \dt{$n$-weakly amenable} if every continuous derivation from $A$ into the $n$th dual space $A^{(n)}$ is inner.
$A$ is \dt{permanently weakly amenable} if it is $n$-weakly amenable for all~$n\in\Nat$\/.

It was shown in \cite{DGG_nWA} that every $C^*$-algebra is permanently weakly amenable, and that every $\LG$ is $n$-weakly amenable for all odd, positive integers~$n$\/.
In \cite{BEJ_free} B.~E.\ Johnson proved that for every free group~$G$\/, the group algebra $\lp{1}(G)$ is $n$-weakly amenable for all even, positive~$n$\/. Combined with \cite[Theorem 4.1]{DGG_nWA}, this shows that for such groups, $\lp{1}(G)$ is permanently weakly amenable.
In an unpublished paper Johnson also showed that for any discrete word-hyperbolic group, the group algebra is permanently weakly amenable.

In fact, for any locally compact group $G$\/, $\LG$ is permanently weakly amenable. Our proof relies heavily on the following result, proved recently by Losert.

\begin{thm}[{\cite[Theorem 1.1]{Los_MG}}]
\label{t:Los_MG}
Let $G$ be a (discrete) group and $X$ a locally compact space on which $G$ has a $2$-sided action by homeo\-morphisms. Then any bounded derivation $D:G\to M(X)$ is inner.
\end{thm}

(The statement in \cite{Los_MG} refers only to those $X$ with a \emph{left} action; however, by standard arguments of Johnson one can reduce the $2$-sided case to the $1$-sided case, see e.g.~\cite[\S2]{BEJ_CIBA}.)

\begin{proof}[Proof that $\LG$ is permanently weakly amenable]
In light of \cite{DGG_nWA} it suffices to show that $\LG$ is $2n$-weakly amenable for all $n\in\Nat$\/.

Let $D:\LG\to \LG^{(2n)}$ be a continuous derivation. By the techniques of \cite[\S1.d]{BEJ_CIBA} $D$ extends to a derivation $\overline{D}:M(G)\to \LG^{(2n)}$\/, where the measure algebra $M(G)$ acts on $\LG^{(2n)}$ through dualizations of its action on $\LG$\/.

Now $\LG^{(2n)}$ is isomorphic, as an $M(G)$-bimodule, to $M(X)$ for some compact space $X$\/.
The action of point masses on $M(X)$ is equivalent to an action of $G$ on $M(X)$\., and $g\mapsto \overline{D}(\delta_g)$ is a bounded derivation from $G$ into $M(X)$\/. Hence by Theorem~\ref{t:Los_MG}\/ this derivation is inner, and this suffices for us to conclude that $\overline{D}:M(G)\to \LG^{(2n)}$ is inner, by \wstar-continuity of $\overline{D}$\/.
\end{proof}

\begin{thm}\label{t:APWA}
Let $G$ be a locally compact SIN group and let $\SG$ be a Segal algebra on $G$\/. Then $\SG$ is approximately permanently weakly amenable (i.e.~for each $n\in\Nat$\/, every continuous derivation $\SG\to \SG^{(n)}$ is approximately inner).
\end{thm}

Note that the case ``$n=0$'' was proved in \cite[Theorem~2.1{$(i)$}]{GhLau02} under the extra hypothesis that our Segal algebra is symmetric.

\begin{proof}
Since $G$ is SIN, it follows from the results of \cite{KR_Seg78} that $\SG$ has a central approximate identity $(e_i)$ which is bounded in the $L^1$-norm.

Let $n\in\Nat$ and let $D:\SG\to\SG^{(n)}$ be a continuous derivation. Our approach is to construct from $D$ a net of continuous derivations $\LG\to\LG^{(n)}$\/, so that we can appeal to Theorem~\ref{t:Los_MG}.

The properties of $(e_i)$\/, together with the derivation property of $D$\/, imply that
\[ D(\SG) \subseteq X_n \defeq \clin\{a\cdot\SG^{(n)}\cdot b \st a,b\in\SG\} \/.\]
Moreover, $(e_i)$ is a two-sided, multiplier-bounded, central approximate identity for $X_n$\/. In particular
\begin{equation}\label{eq:SIDESTEP}
 \lim_i e_i^2\cdot D(f) =D(f) \qquad(f\in\SG)
\end{equation}
where the limit is taken in the norm topology of $\SG^{(n)}$\/.

For each $i$\/, define a continuous linear mapping $\tau_i:\LG\to \SG$ by
\[ \tau_i(f)=f*e_i \qquad(f\in \LG) \]
and let $\theta:\SG\to\LG$ denote the (continuous) inclusion map. Both $\tau_i$ and $\theta$ are left 
 $\LG$-module morphisms and are also $\SG$-bimodule morphisms. Clearly, $\tau_i\theta(f)=f*e_i$ for $f\in\SG$\/; so by induction, for each $n\in\Nat$ the map $(\tau_i\theta)^{(n)}:\SG^{(n)}\to\SG^{(n)}$ satisfies
\[ (\tau_i\theta)^{(n)}(F) = \left\{ \begin{aligned}
F\cdot e_i & \quad\text{if $n$ is even} \\
e_i\cdot F & \quad\text{if $n$ is odd} 
\end{aligned} \right.\qquad(F\in\SG).
\]

Define $\Delta_i: L^1(G)\to L^1(G)^{(n)}$ by
\[ \Delta_i(f)=\left\{\begin{aligned}
\theta^{(n)}[ D(f*e_i) - f\cdot D(e_i)] & \quad\text{if $n$ is even} \\ 
\tau_i^{(n)}[ D(e_i*f) - D(e_i)\cdot f] & \quad\text{if $n$ is odd}
\end{aligned}\right.\qquad(f\in L^1(G)).
\]
Then $\Delta_i$ is a continuous linear map, and for $f\in\SG$ we have
\begin{equation}\label{eq:Branch3}
\Delta_i(f) = \left\{\begin{aligned}
\theta^{(n)}(D(f)\cdot e_i) & \quad\text{if $n$ is even} \\
\tau_i^{(n)}(e_i\cdot D(f)) & \quad\text{if $n$ is odd.}
\end{aligned}\right.
\end{equation}
Since $e_i$ is central, it is straightforward to verify using \eqref{eq:Branch3} that $\Delta_i(f*g)=f\cdot\Delta_i(g) +\Delta_i(f)\cdot g$ for all $f,g\in\SG$\/. Therefore, since $\SG$ is dense in $\LG$\/, it follows that $\Delta_i$ is a derivation from $\LG$ to $\LG^{(n)}$\/.

By Theorem \ref{t:Los_MG} there exists $m_i\in \LG^{(n)}$ such that
\[ \Delta_i(f)=f\cdot m_i -m_i\cdot f \qquad(f\in\LG).\]
In particular, for $f\in\SG$\/, Equation~\eqref{eq:Branch3} implies that for even $n$ we have
\begin{align*}
D(f)\cdot e_i^2
 & = (\tau_i\theta)^{(n)}(D(f)\cdot e_i) \\
 & = \tau_i^{(n)}\Delta_i(f) \\
 & = f\cdot\tau_i^{(n)}(m_i) - \tau_i^{(n)}(m_i)\cdot f\/,
\end{align*}
while for odd $n$ we have
\begin{align*}
e_i^2\cdot D(f)
 & = (\tau_i\theta)^{(n)} (e_i\cdot D(f)) \\
 & = \theta^{(n)}\Delta_i(f) \\
 & = f\cdot\theta^{(n)}(m_i) -\theta^{(n)}(m_i)\cdot f\,.
\end{align*}

Take $n_i = \tau_i^{(n)}(m_i)$ if $n$ is even, and take $n_i=\theta^{(n)}(m_i)$ if $n$ is odd.
Then $n_i\in\SG^{(n)}$ for all~$i$\/, and for every $f\in\SG$ we have, by \eqref{eq:SIDESTEP},
\[ D(f)=\lim_i e_i^2\cdot D(f) =\lim_i f\cdot n_i - n_i \cdot f\]
Thus $D$ is approximately inner, as required.
\end{proof}

\begin{remstar}
Our construction actually provides a \emph{bounded} net of inner derivations which approximate $D$\/, although the net of implementing elements need not be bounded.
\end{remstar}

\end{section}

\begin{section}{$\lp{1}$-convolution algebras of totally ordered sets}\label{s:l1ordered}
Recall that a \dt{semilattice} is a commutative semigroup in which every element is idempotent. The $\lp{1}$-convolution algebras of semilattices provide interesting examples of commutative Banach algebras.
However, amenability is too strong a notion for such algebras: if $S$ is a semilattice then the convolution algebra $\Alg{S}$ is amenable if and only if $S$ is \emph{finite}~\cite[Theorem~10]{DuncNam}.

It is not clear to the authors exactly which semilattices have \appam\ $\lp{1}$-convolution algebras. In the case where the semilattice is \emph{totally ordered} we can do better.

Let $\Lm$ be a non-empty, totally ordered set, and regard it as a semigroup by defining the product of two elements to be their maximum. The resulting semigroup, which we denote by $\Lmmax$, is a semilattice.
We may then form the $\lp{1}$-convolution algebra $\Alg{\Lmmax}$. For every $t\in\Lmmax$ we denote the point mass concentrated at $t$ by~$e_t$\/. The definition of multiplication in $\Alg{\Lmmax}$ ensures that $e_se_t=e_{\max(s,t)}$ for all $s$ and~$t$\/.

\begin{remstar}
One could also turn $\Lm$ into a semilattice $\Lm_\wedge$ by defining the product of two elements to be their \emph{minimum}. This is in some sense more natural, for reasons we shall not discuss here; we have chosen to work with $\Lmmax$ as this fits our main example (in Theorem \ref{t:total-ord_not_CAC}) better.
\end{remstar}

\begin{thm}\label{t:total-ord_is_BAC}
Let $\Ind$ be \emph{any} totally ordered set. Then $\Alg{\Imax}$ is boundedly \appctr.
\end{thm}

\begin{remstar}
The special case of $\Ind=\Nat$ or $\Nat^{\op}$ was done in \cite{GhLZ_genam2}. Our arguments are a more abstract version of the ones there.
\end{remstar}

We prove the theorem in several steps. First, by following the proof of \cite[Theorem~5.10]{GhLZ_genam2}, it suffices to prove that $\fu{\Alg{\Imax}}$ has a multiplier-bounded approximate diagonal, in the sense of Definition~\ref{dfn:MBAD}. Moreover, we can identify $\fu{\Alg{\Imax}}$ with $\Alg{\wtild{\Ind}_\vee}$\/, where $\wtild{\Ind}$ denotes the disjoint union of $\Ind$ with an adjoined least element. Clearly $\wtild{\Ind}$ is also a totally ordered set, and so to prove Theorem \ref{t:total-ord_is_BAC} it suffices to prove that
\begin{quote}
{\it for any totally ordered set $\Ind$, $\Alg{\Imax}$ has a multiplier-bounded approximate diagonal.}
\end{quote}

It is useful to first consider the case of a finite totally ordered set. More precisely, let $F$ be a finite subset of~$\Ind$, and enumerate its elements in increasing order as
\[ \min(F) = c(0)<c(1)<\dots< c(n) =\max(F)\]
say. We then define $\Delta_F\in \Alg{\Imax}\ptp\Alg{\Imax}$ by
\begin{equation}\label{eq:dfn_findiag}
\Delta_F = \left(\sum_{j=1}^n (e_{c(j-1)}-e_{c(j)})\tp(e_{c(j-1)}-e_{c(j)})
\right) + e_{c(n)}\tp e_{c(n)}\/.
\end{equation}
A small calculation shows that $\pi(\Delta_F) = e_{c(0)}$\/, so that
\begin{equation}\label{eq:findiag1}
e_\lm\pi(\Delta_F)=e_\lm \quad\text{ for all $\lm\in F$.}
\end{equation}
It is also easily checked that
\begin{equation}\label{eq:findiag2}
e_\lm\cdot \Delta_F = \Delta_F\cdot e_\lm \quad\text{ for all $\lm\in F$.}
\end{equation}
and thus $\Delta_F$ is a diagonal for the subalgebra $\Alg{F_\vee}\subseteq \Alg{\Imax}$.

Having seen how to construct a diagonal for the finite case, we now proceed to the general case. Let $\Fin$ be the set of all non-empty \emph{finite} subsets of
$\Ind$, and order $\Fin$ with respect to inclusion, so that for
any $E$ and $F$ in $\Fin$, $E\preceq F$ if and only if $E\subseteq F$\/.

The following result will, by the remarks above, imply Theorem \ref{t:total-ord_is_BAC}.
\begin{propn}\label{p:MBAD-for-Imax}
\label{p:D}
The net $(\Delta_F)_{F\in \Fin}$ is a multiplier-bounded approximate diagonal for $\Alg{\Imax}$.
\end{propn}

We isolate the key technical estimate as a lemma.
\begin{lemma}\label{l:Marsha}
Let $b\in \Alg{\Imax}$, $F\in\Fin$\/. Then $\norm{b\cdot\Delta_F-\Delta_F\cdot b}\leq 6\norm{b}$\/.
\end{lemma}

\begin{proof}
By the triangle inequality and the definition of the $\lp{1}$-norm, we can without loss of generality assume that $b=e_\lm$ for some $\lm\in\Ind$\/. Thus it suffices to prove that
\begin{equation}\label{eq:MBAD1}
 \norm{e_\lm\cdot \Delta_F-\Delta_F\cdot e_\lm}\leq 6 \quad\text{ for all $F\in\Fin$\/.}
\end{equation}
This estimate holds trivially if $F$ consists of only one point, so we shall henceforth assume that $\abs{F}\geq 2$\/.

As before we enumerate the elements of $F$ in increasing order as $c(0)<c(1)<\ldots<c(n)$\/.
We consider three possibilities. If $\lm \geq c(n)$\/, then $e_\lm\cdot\Delta_F=e_\lm\tp e_{c(n)}$ and $\Delta_F\cdot e_\lm= e_{c(n)}\tp e_\lm$, so that
\[ \norm{e_\lm\cdot\Delta_F-\Delta_F\cdot e_\lm} \leq 2\/, \]
and so \eqref{eq:MBAD1} certainly holds. At the other extreme, if $\lm\leq c(0)$ then $e_\lm\cdot\Delta_F= \Delta_F=\Delta_F\cdot e_\lm$, so that \eqref{eq:MBAD1} once again holds.

The third remaining possibility is that $c(0) < \lm < c(n)$\/. Let 
\[ m=\min\{k \st c(k) > \lm\}\]
so that $1\leq m\leq n$ and $c(m-1) < \lm <c(m)$\/.

When we calculate $e_\lm\cdot\Delta_F-\Delta_F\cdot e_\lm$ using the formula~\eqref{eq:dfn_findiag}, most of the terms cancel and we obtain
\[\begin{aligned}
 e_\lm\cdot\Delta_F - \Delta_F\cdot e_\lm
 & = \left\{\begin{gathered}
	e_\lm(e_{c(m-1)}-e_{c(m)})\tp(e_{c(m-1)}-e_{c(m)}) \\
	- (e_{c-1(m)}-e_{c(m)})\tp(e_{c(m-1)}-e_{c(m)})e_\lm
	\end{gathered}\right. \\
 & = \left\{\begin{gathered}
	(e_\lm -e_{c(m)})\tp (e_{c(m-1)}-e_{c(m)}) \\
	- (e_{c(m-1)}-e_{c(m)})\tp (e_\lm - e_{c(m)})\/.
	\end{gathered}\right.
\end{aligned} \]
Expanding out and using the triangle inequality gives
$\norm{e_\lm\cdot\Delta_F-\Delta_F\cdot e_\lm} \leq 6$ as required.
\end{proof}

\begin{proof}[Proof of Proposition \ref{p:MBAD-for-Imax}]
Fix $a\in\Alg{\Imax}$.
We have already seen in Lemma \ref{l:Marsha} that $\norm{a\cdot\Delta_F-\Delta_F\cdot a}\leq 6\norm{a}$ for every $F\in\Fin$\/. Also, since $\pi(\Delta_F)=e_{\min(F)}$, we have
\begin{equation}\label{eq:Amber}
\norm{a\pi(\Delta_F)-a} \leq 2\norm{a} \quad\text{ for every $F$ in $\Fin$\/.}
\end{equation}
Thus the `multiplier-bounded' part of the defining condition \eqref{eq:MBAD-dfn} is satisfied.

It remains to show that, given $\veps >0$, there exists $F_0\in\Fin$ such that\[ \norm{a\pi(\Delta_F)-a} <\veps \text{ and } \norm{a\cdot\Delta_F-\Delta_F\cdot a} < \veps\]
for any $F\in\Fin$ with $F\supseteq F_0$.

Fix $\veps>0$ and choose $F_0\in\Fin$ such that $\sum_{\lm \in\Ind\setminus F_0} \abs{a_\lm} \leq \veps/6$, and let $F\in\Fin$ with $F\supseteq F_0$. Let $\tilde{a}$ denote the obvious truncation of $a$ to the subset $F$
(i.e.~$\tilde{a}_\lm=a_\lm$ if $\lm\in F$ and
$\tilde{a}_\lm=0$ otherwise). Note that
$\norm{a-\tilde{a}}\leq\veps/6$.

Since $\tilde{a}\in\Alg{F_\vee}$, we deduce from Equation~\eqref{eq:findiag2} and the estimate \eqref{eq:Amber} that
\[\norm{a\cdot\Delta_F-\Delta_F\cdot a} = \norm{(a-\tilde{a})\cdot\Delta_F-\Delta_F\cdot(a-\tilde{a})} \leq 6\norm{a-\tilde{a}}\leq\veps \]
Finally, using Equation~\eqref{eq:findiag1} and Lemma \ref{l:Marsha} we obtain
\[ \norm{a\pi(\Delta_F)-a}
	 =\norm{(a-\tilde{a})\pi(\Delta_F)-(a-\tilde{a})} 
	 \leq 2\norm{a-\tilde{a}} =\veps/3 \/, \]
and the proof is complete.
\end{proof}

\begin{remstar}
If the set $\Ind$ is countable, then the net $(\Delta_F)_{F\in\Fin}$ has a subnet which is a sequence (take any enumeration of $\Ind$ as $\{t_1,t_2,\ldots\}$ and let $\wtild{\Delta}_n\defeq \Delta_{\{t_1,\ldots,t_n\}}$\/). So if $\Ind$ is countable then $\Alg{\Imax}$ is sequentially \appctr. 
\end{remstar}

While sequential \appamy\ implies bounded \appamy, the converse is false. This is proved by combining Theorem \ref{t:total-ord_is_BAC} with the following result.
\begin{thm}\label{t:total-ord_not_CAC}
Let $\Lm$ be an uncountable well-ordered set. Then $\Alg{\Lmmax}$ is not
 sequentially \appam.
\end{thm}

(Recall that a totally ordered set is \dt{well-ordered} if every non-empty subset has a least element: well-ordered sets are precisely those ordered sets which are order-isomorphic to ordinals.)

In proving Theorem \ref{t:total-ord_not_CAC} we shall use some basic facts on the character theory of $\Alg{\Lmmax}$.
It is clear that the characters on $\Alg{\Lmmax}$ correspond to the non-zero semigroup homomorphisms from $\Lmmax$ to the two-element semigroup $\{0,1\}$; and a little thought gives the following characterization.

\begin{propn}
The characters on $\Alg{\Lmmax}$ are all of the form $\indic{\Lm\setminus U}$, where $U$ is a proper (and possibly empty) subset of $\Lm$ that is \emph{upwards-directed} with respect to the given order on $\Lm$.
\end{propn}

\begin{eg}
Take $\Lm$ to be the real line with its usual ordering. Then the characters on $\Alg{\Lmmax}$ are either of the form $\indic{(-\infty, t)}$ or $\indic{(-\infty, t]}$\/.
\end{eg}


If $U$ is a non-empty, upwards-directed subset of a well-ordered set $\Lm$, then $U$ has a least element, $u$ say: hence $U=\{x\in\Lm : x\geq u\}$. Thus the complements of upwards-directed sets are all of the form $\{y: y< u \}$\/.

If $\lm$ is an element of a well-ordered set and it is not maximal, then there is a unique minimal element greater than $\lm$, which we shall denote by $\lm+1$.

\begin{notn}
Let $\Lm$ be a well-ordered set and consider the algebra $\Alg{\Lmmax}$. If $\lm\in \Lm$ we denote by $\wtild{\lm}$ the character $\indic{< \lm}$. If $\lm$ is maximal in $\Lm$ then we adopt the convention that $\wtild{\lm+1}$ is the \dt{augmentation character} $\indic{\Lm}$.
\end{notn}

The following is then obvious: we isolate it as a lemma for later reference.
\begin{lemma}\label{l:diff-char}
Let $\Lm$ be a well-ordered set and let $\lm\in\Lm$. Then \[ \delta_\lm = \wtild{\lm+1}-\wtild{\lm}\/,\]
where $\delta_\lm$ denotes the point mass at $\lm$, regarded as an element of $\Alg{\Lmmax}^*$.
\end{lemma}

Our proof of Theorem \ref{t:total-ord_not_CAC} uses our earlier observations on the characters of $\Alg{\Lmmax}$, together with Lemma \ref{l:Gelf-tr}. Intuitively, the idea is that the Gelfand transforms of elements in $\Alg{\Lmmax}\ptp\Alg{\Lmmax}$ are bad approximations to the indicator function of the set $\{(\lm,\lm)\st \lm\in\Lm\}$, so that if $\Lm$ is uncountable then no countable net $(\Delta_n)$ can have the properties described in Lemma \ref{l:Gelf-tr}. We make this idea precise as follows. 

\begin{lemma}\label{l:C}
Let $\Ind$ be an uncountable index set and let $(F_n)$ be a countable family in $\lp{1}(\Ind\times\Ind)^{**}$. Then there exist uncountably many $t\in\Ind$ such that
\[ \pair{F_n}{\delta_t\tp\delta_t}=0 \quad\hbox{ for all $n$}. \]
\end{lemma}
\begin{proof}
In view of the direct sum decomposition
\[ \lp{1}(\Ind\times\Ind)^{**} = \lp{1}(\Ind\times\Ind)\oplus c_0(\Ind\times\Ind)^{\perp} \]
we may write each $F_n$ as $\kp(f_n)+G_n$ where $G_n\in c_0(\Ind\times\Ind)^\perp$, $f_n\in\lp{1}(\Ind\times\Ind)$ and $\kp$ is the natural embedding of $\lp{1}(\Ind\times\Ind)$ in its bidual.

Let
\[ S=\bigcup_n \{ t\in\Ind \st f_n(t,t)\neq 0 \}  \]
Since each $f_n$ has \emph{countable support}, $S$ is countable. In particular
$\Ind\setminus S$ is uncountable, and for any $t\in \Ind\setminus S$ we have
\[ \pair{F_n}{\delta_t\tp\delta_t} = (f_n)_{t,t} + \pair{G_n}{\delta_t\tp\delta_t} =0\]
as claimed.
\end{proof}

\begin{proof}[Proof of Theorem \ref{t:total-ord_not_CAC}] Suppose $\Alg{\Lmmax}$ is sequentially \appam. Since $\Lm$ is an ordinal it has a least element, and consequently $\Alg{\Lmmax}$ has an identity element. Hence Lemma~\ref{l:Gelf-tr} applies and there is a sequence $\Delta_n\in \left(\Alg{\Lmmax}\ptp\Alg{\Lmmax}\right)^{**}$ such that
\begin{equation}\label{eq:SHAZAM}
 \pair{\Delta_n}{\varphi\tp\varphi} = 1 \quad\hbox{ for all $n$ and }\quad \lim_n\pair{\Delta_n}{\varphi\tp\chi} = 0
\end{equation}
for every pair of distinct characters $\varphi,\chi$.

By Lemma \ref{l:C} there exists $\lm\in\Lambda$ such that 
 \[ \pair{\Delta_n}{\delta_\lm\tp\delta_\lm} =0 \quad\text{ for all $n$} \]
and hence by Lemma \ref{l:diff-char}
\[ 0 = \left\{ \begin{gathered}
	\pair{\Delta_n}{\wtild{\lm}\tp\wtild{\lm}} -
	\pair{\Delta_n}{\wtild{\lm}\tp\wtild{\lm+1}} \\
	 - \pair{\Delta_n}{\wtild{\lm+1}\tp\wtild{\lm}} + \pair{\Delta_n}{\wtild{\lm+1}\tp\wtild{\lm+1}}
	\end{gathered}\right. \]
But by Equation \eqref{eq:SHAZAM} the right-hand side converges to $2$ as $n\to\infty$, which is a flagrant contradiction.
\end{proof}

\begin{remstar}
The proof just given yields something formally stronger, namely that $\Alg{\Lmmax}$ cannot have an \appdiag\ with \emph{countable} indexing set. We do not pursue this further in this paper, chiefly because we know of no Banach algebra which has a countably-indexed \appdiag\ and yet has no sequential \appdiag.
\end{remstar}

\end{section}
\begin{section}{Algebras of pseudo-functions on discrete groups}\label{s:gpalg}

Let $\Gm$ be a discrete group, with convolution algebra $\lp{1}(\Gm)$. Given $p\in(1,\infty)$ we may consider the left regular representation of $\Gm$ on $\lp{p}(\Gm)$, and this gives an injective continuous algebra homomorphism
\[ \theta_p:\lp{1}(\Gm)\to \Bdd(\lp{p}(\Gm))\/. \]
We denote by $\psf{p}(\Gm)$ the norm-closure in $\Bdd(\lp{p}(\Gm))$ of the range of $\theta_p$. Note that $\psf{2}(\Gm)$ is nothing but the reduced $C^*$-algebra of $\Gm$.

 If $\Gm$ is amenable, then by Johnson's theorem the convolution algebra $\lp{1}(\Gm)$ is amenable, and since amenability is inherited by closures under Banach algebra norms we deduce that $\psf{p}(\Gm)$ is amenable; in particular the reduced $C^*$-algebra $\redC{\Gm}$ is amenable. The converse result -- that amenability of $\redC{\Gm}$ implies amenability of~$\Gm$ -- was proved by Bunce\footnotemark\ in~\cite{Bun_amenCG}.
\footnotetext{It had already been shown by Lance that \emph{nuclearity} of $\redC{\Gm}$ implies amenability of~$\Gm$.}
With some modifications one can adapt his proof to show that amenability of any one of the algebras $\psf{p}(\Gm)$ is enough to force amenability of $\Gm$.

In \cite{GhL_genam1} it was shown that \appamy\ of the group algebra $L^1(G)$ implies amenability of $G$, by generalizing the well-known argument for \emph{amenability} of $\Lp{1}(G)$. We shall now show that by combining arguments from \cite{Bun_amenCG} and \cite{GhL_genam1} we have the following theorem.

\begin{thm}\label{t:AA_PF(G)}
Let $\Gm$ be a discrete group. Then the following are equivalent:
\begin{itemize}
\item[$(i)$] $\Gm$ is amenable;
\item[$(ii)$] $\psf{p}(\Gm)$ is amenable for all $p\in (1,\infty)$;
\item[$(iii)$] $\psf{p}(\Gm)$ is \appam\ for some $p\in (1,\infty)$.
\item[$(iv)$] $\psf{p}(\Gm)$ is pseudo-amenable for some $p\in (1,\infty)$.
\end{itemize}
\end{thm}

As mentioned above, the implications $(i)\iff(ii)$ are already known, while the implication $(ii)\implies (iii)$ is trivial; the implication $(iii)\implies(iv)$ follows from \cite[Proposition 3.2]{GhZ_pseudo} and only uses the fact that $\psf{p}(\Gm)$ has an identity element.
 Therefore our contribution is to prove the implication $(iv)\implies (i)$. Taking $p=2$, our proof will give a slightly streamlined version of Bunce's arguments, in that we are able to forgo technical arguments with states on $C^*$-algebras in favour of more basic positivity arguments with measures on compact spaces: that is, we can make do with commutative rather than noncommutative measure theory.

Our idea is to follow Bunce's construction up to the point where he produces, from the assumption that $\redC{\Gm}$ is \emph{amenable}, a non-zero element $\rho$ in $\lp{\infty}(\Gm)^*$ which satisfies\footnotemark
\footnotetext{In \cite{Bun_amenCG} $\rho$ is described as being `left-invariant'\/: we adopt the opposite and more usual convention, and say $\rho$ is \emph{right-invariant}.}
\[ \rho(g\cdot f) = \rho(f)\qquad\text{ for all $f\in\lp{\infty}(\Gm)$\/.}\]
In our setting we merely obtain a net $(\phi_\al)$ of functionals on $\lp{\infty}(\Gm)$ which satisfies
\[ \phi_\al({\mathbf 1})\to 1 \quad\text{ and } \quad
  \norm{\phi_\al\cdot g -\phi_\al} \to 0 \text{ for each $g\in\Gm$.}\]
We then use this net to obtain a ``genuine'' invariant mean on $\lp{\infty}(\Gm)$, by following the last part of the proof of \cite[Theorem 3.2]{GhL_genam1}. For convenience we isolate the relevant argument and state it as the following lemma.

\begin{lemma}\label{l:GhL_Ringrose}
Let $G$ be a locally compact group, and let $T$ be a compact $G$-space on which $G$ acts from the right by homeomorphisms. Equip $M(T)$ with its usual norm, and regard it as a right Banach $G$-module.

Suppose we have a net $(\varphi_i)$ of Radon measures on $T$ which satisfies the following conditions:
\begin{enumerate}[$(i)$]
\item $\inf_i \norm{\varphi_i} >0$;
\item $\norm{\varphi_i\cdot g - \varphi_i} \to 0$ for all $g\in G$.
\end{enumerate}
Then there exists a probability measure $n$ on $T$ such that $n\cdot g=n$.
\end{lemma}

We give the proof for sake of completeness (cf.~the proof of \cite[Theorem 3.2]{GhL_genam1}).
\begin{proof}
Set $n_i = \norm{\varphi_i}^{-1}\varphi_i$. The hypothesis that $\norm{\varphi_i}$ is bounded below then ensures that
\[ \norm{n_i\cdot g - n_i} \to 0 \qquad\text{ for all $g\in G$}. \]

For any two Radon measures $\mu$, $\nu$ on $T$ we have
$\norm{\mu-\nu} \geq \norm{ \abs{\mu}-\abs{\nu} }$,
an inequality which can easily be deduced from the definition of the total variation of a measure.
Therefore, since $\abs{\mu\cdot g}=\abs{\mu}\cdot g$ for any $\mu\in M(T)$, we have
\[ \norm{\abs{n_i}\cdot g -\abs{n_i} } =\norm{ \abs{n_i\cdot g} -\abs{n_i}} \leq \norm{n_i\cdot g - n_i} \to 0 \,, \]
for every $g\in G$.

Take $n$ to be any $w^*$-cluster point of the net $(\abs{n_i})$\/. Since $\abs{n_i}({\mathbf 1})=1$ for all $i$\/, we have $n({\mathbf 1}) = 1$\/; 
and for any $g\in G$ and $f\in C(T)$\/, we have
\[ \abs{(n\cdot g - n)(f)} \leq \limsup_i \abs{(\abs{n_i}\cdot g-\abs{n_i})(f)} =0 \,\]
so that $n\cdot g =n$ for all $g\in G$\/.
\end{proof}

\begin{proof}[Proof of Theorem \ref{t:AA_PF(G)}, $(iv)\implies (i)$]
Our aim is to construct a right-invariant mean on $\lp{\infty}(\Gm)$\/. To fix notation, we recall that the usual left action of $\Gm$ on $\lp{\infty}(\Gm)=\lp{1}(\Gm)^*$ is defined by
\[ (g\cdot f)(x) = f(g^{-1}x)\qquad\text{for $f\in\lp{\infty}(\Gm)$ and $g, x\in \Gm$\/.} \] 
For each $g\in\Gm$ let $L_g$ be the isometric, invertible operator on $\lp{p}(\Gm)$ given by left translation,~i.e.
\[ (L_g k)(x) = k(g^{-1}x)\quad\text{for all $k\in\lp{p}(\Gm)$.}\]

We regard $\psf{p}(\Gm)$ as a subalgebra of $B(\lp{p}(\Gm))$. Take $\tau$ to be the functional given by
\[ \tau(T) = \pair{T\delta_e}{\delta_e} \quad\quad(T\in \Bdd{\lp{p}(\Gm)})\]
where $\delta_e$ is the basis vector of $\lp{p}(\Gm)$ that takes the value $1$ at $e$ and the value $0$ everywhere else. Clearly $\tau(I)=1$.
A simple computation shows that for any $a,b$ in the group algebra $\Cplx\Gm$, we have
\[ \tau(\theta_p(a)\theta_p(b)) = \tau(\theta_p(b)\theta_p(a)) \, \]
and so by continuity the restriction of $\tau$ to $\psf{p}(\Gm)$ defines a non-zero trace.

Suppose that $\psf{p}(\Gm)$ is pseudo-amenable. By Lemma~\ref{l:extension}, there exists a net $(\psi_\al)$ in $B(\lp{p}(\Gm))$ such that
\begin{equation}\label{eq:nonzero}
\lim_\al \psi_\al(I)=1
\end{equation}
 and
\[ \lim_\al \sup_{T\in B(\lp{p}(\Gm)), \norm{T}\leq 1} \abs{\psi_\al(aT-Ta)}=0 \quad\text{ for all $a\in \psf{p}(\Gm)$.}\]
In particular, for any $g\in \Gm$ we have
\begin{equation}\label{eq:conj-inv}
\begin{gathered}
 \sup_{M\in B(\lp{p}(\Gm)),\norm{M}\leq 1} \abs{\psi_\al(L_gML_{g^{-1}}) - \psi_\al(M)} \\
 \leq \sup_{T\in B(\lp{p}(\Gm)), \norm{T}\leq 1} \abs{\psi_\al(L_gT-TL_g)} \to 0
\end{gathered}
\end{equation}

Regard $\lp{\infty}(\Gm)$ as an algebra with pointwise multiplication and supremum norm. There is an embedding of $\lp{\infty}(\Gm)$ as a closed unital subalgebra of $B(\lp{p}(\Gm))$, defined by sending a bounded function $f\in\lp{\infty}(\Gm)$ to the ``diagonal multiplication'' operator $M_f$ where
\[ (M_f k)(x) = f(x)k(x) \qquad(k\in\lp{p}(\Gm), x\in\Gm).\]
Then a direct calculation shows that
\begin{equation}\label{eq:crossprod}
 M_{(g\cdot f)} = L_g M_f (L_g)^{-1} \qquad\text{ for all $f\in\lp{\infty}(\Gm)$ and $g\in \Gm$} \/.
 \end{equation}

For each $\al$ define $\phi_\al\in \lp{\infty}(\Gm)^*$ by $\phi_\al(f) = \psi_\al(M_f)$, $f\in \lp{\infty}(\Gm)$\/. It follows from Equations \eqref{eq:nonzero}, \eqref{eq:conj-inv} and \eqref{eq:crossprod} that
\[ \lim_\al \phi_\al({\mathbf 1})=1 \]
and
\[ \lim_\al \norm{\phi_\al\cdot g -\phi_\al} = \sup_{f\in\lp{\infty}(\Gm),\norm{f}\leq1} \abs{\phi_\al(g\cdot f)-\phi_\al(f)} =0 \quad\text{ for all $g\in \Gm$.}\]

To finish we observe that $\lp{\infty}(\Gm)$ may be identified with the space of continuous functions on a \emph{compact} $\Gm$-space $T$ (namely, take $T$ to be the Stone-\v{C}ech compactification of $\Gm$), and hence we may identify each $\phi_\al$ with a Radon measure on $T$. By Lemma~\ref{l:GhL_Ringrose}, there exists a positive functional $n\in\lp{\infty}(\Gm)^*$ satisfying $n({\mathbf 1})=1$ and $n\cdot g=n$ for all $g\in \Gm$, and hence $\Gm$ is amenable as claimed.
\end{proof}

Specializing to the case $p=2$ (i.e.~the reduced $C^*$-algebra $\redC{\Gm}$), we have the following corollary.

\begin{coroll}
The full group $C^*$-algebra $C^*(\Gm)$ is \appam\ if and only if $\Gm$ is amenable.
\end{coroll}
\begin{proof}
We first recall without proof some basic facts about $C^*(\Gm)$: firstly, it is by definition the completion of $\lp{1}(\Gm)$ in a certain $C^*$-norm; and secondly, there is a canonical quotient homomorphism from $C^*(\Gm)$ onto $\redC{\Gm}$.

Now, suppose $\Gm$ is amenable: then $\lp{1}(\Gm)$ is amenable. As just mentioned, the inclusion homomorphism $\lp{1}(\Gm)\to C^*(\Gm)$ is continuous with dense range, and therefore $C^*(\Gm)$ must also be amenable.

Conversely, suppose that $C^*(\Gm)$ is \appam. By \cite[Propos\-ition~2.2]{GhL_genam1},
 \appamy\ passes to quotient algebras, and so $\redC{\Gm}$ is \appam. Now apply Theorem \ref{t:AA_PF(G)} in the case $p=2$.
\end{proof}

\begin{remstar}
Using the fact that the canonical tracial state $\tau$ on $\redC{\Gm}$ actually extends to a tracial state on the von Neumann algebra $VN(\Gm)$, we can adapt the proof of Theorem \ref{t:AA_PF(G)} to show the following result:
\begin{quote}
{\it If $A$ is a closed unital subalgebra of $VN(\Gm)$, with $\redC{\Gm}\subseteq A$, and furthermore $A$ is \appam, then $\Gm$ is amenable.}
\end{quote}
\end{remstar}

\end{section}

\appendix

\begin{section}{Implications}

\subsection{General implications}
\[ \begin{diagram}[nohug]
& & & & \text{biflat} \\
& & & & \uTo \\
AA & \lTo^{\dfn} & BAA & \lTo^{\dfn} & \text{amenable} & \rTo^{\dfn} & PsA \\
\uTo^{\dfn}\dTo_{\rm\cite[Thm 2.1]{GhLZ_genam2}} & & \uTo^{\dfn} & \ldTo^{\rm \cite{FMG_approx}} & \uTo_{\dfn} &  & \uTo^{\dfn} \dNo_{(\star)} \\
AC & \lTo_{\dfn} & BAC & \lTo_{\dfn} & \fbox{contractible} & \rTo_{\dfn} & PsC
\end{diagram} \]
Here the counterexample $(\star)$ to `pseudo-amenable implies pseudo-contractible' follows because unital pseudo-contractible algebras must be contractible \cite[Theorem~2.4]{GhZ_pseudo}, while there are unital pseudo-amenable algebras which are not even amenable: perhaps the simplest example is $\lp{1}(\Nat,\max)$.

\subsection{Commutative settings}
\[ \begin{diagram}[nohug]
CAI+AA & \rTo^{\rm\cite[Prop~3.3]{GhZ_pseudo}} & PsA & \lTo^{\dfn} & \Cmm+PsA \\
\uTo^{\rm\cite[Lem~2.2]{GhL_genam1}} & & & & \uTo_{\dfn} \\
\Cmm+AA & & \lNo_{\rm\cite[Thm~4.1]{DLZ_St06}}  & & \fbox{$\Cmm$+ PsC}
\end{diagram} \]
Here the counterexample to `pseudo-contractible implies \appam' is given by $\lp{1}(\Nat)$ with pointwise multiplication.

\subsection{Approximate identities}
\begin{diagram}[nohug]
PsA+BAI & \ltwoway^{\rm \cite[Prop~3.2]{GhZ_pseudo}} & AA+BAI & \lTo^{\dfn}   & BAA+MBAI \\
 & & \uTo & &  \uTo^{\dfn}\dTo_{\bf\ref{t:BAA+MBAI}} \\
BAC & \pile{\lTo^{\dfn}\\ \rTo_{\bf\ref{c:BAC-has-BAI}}} & \fbox{BAC+BAI} & \rTo^{\dfn} & BAA+BAI 
\end{diagram}

\end{section}

\vfill\eject

\newcommand{\affil}[1]{\vspace{1.0em}\parbox{\textwidth}{#1}}
\newcommand{\email}[1]{{\tt #1}}

\affil{\sl
Y. Choi, F. Ghahramani and Y. Zhang \\
Department of Mathematics\\
University of Manitoba\\
Winnipeg, Canada R3T 2N2
}

\medskip

\email{y.choi.97@cantab.net}

\email{fereidou@cc.umanitoba.ca}

\email{zhangy@cc.umanitoba.ca}

\end{document}